\documentclass[12pt]{amsart}
\usepackage{amscd}
\def\frk{\frak}               

\def\Phi{{\frk n}}
\def\Phi{{\frk N}}
\def\opn#1#2{\def#1{\operatorname{#2}}} 
\opn\chara{char}
\opn\length{\ell}
\opn\pd{pd}
\opn\rk{rk}
\opn\projdim{proj\,dim}
\opn\injdim{inj\,dim}
\opn\rank{rank}
\opn\depth{depth}
\opn\grade{grade}
\opn\height{height}
\opn\embdim{emb\,dim}
\opn\codim{codim}

\opn\Tr{Tr}
\opn\bigrank{big\,rank}
\opn\superheight{superheight}\opn\lcm{lcm}
\opn\trdeg{tr\,deg}%
\opn\reg{reg}
\opn\lreg{lreg}
\opn\ini{in}
\opn\div{div}
\opn\Div{Div}
\opn\cl{cl}
\opn\Cl{Cl}
\opn\Spec{Spec}
\opn\Supp{Supp}
\opn\supp{supp}
\opn\Sing{Sing}
\opn\Ass{Ass}
\opn\Ann{Ann}
\opn\Rad{Rad}
\opn\Soc{Soc}
\opn\Ker{Ker}
\opn\Coker{Coker}
\opn\Am{Am}
\opn\Hom{Hom}
\opn\Tor{Tor}
\opn\Ext{Ext}
\opn\End{End}
\opn\Aut{Aut}
\opn\id{id}

\opn\nat{nat}
\opn\pff{pf}
\opn\Pf{Pf}
\opn\GL{GL}
\opn\SL{SL}
\opn\mod{mod}
\opn\ord{ord}
\opn\Gin{Gin}
\opn\aff{aff}
\opn\con{conv}
\opn\relint{relint}
\opn\st{st}
\opn\lk{lk}
\opn\cn{cn}
\opn\core{core}
\opn\vol{vol}
\opn\link{link}
\opn\star{star}
\opn\sort{sort}
\opn\gr{gr}

\def\pot#1#2{#1[\kern-0.28ex[#2]\kern-0.28ex]}

\opn\dirlim{\underrightarrow{\lim}}
\opn\inivlim{\underleftarrow{\lim}}
\def\Implies{\ifmmode\Longrightarrow \else
     \unskip${}\Longrightarrow{}$\ignorespaces\fi}
\def\implies{\ifmmode\Rightarrow \else
     \unskip${}\Rightarrow{}$\ignorespaces\fi}
\def\iff{\ifmmode\Longleftrightarrow \else
     \unskip${}\Longleftrightarrow{}$\ignorespaces\fi}

\let\:=\colon
\newtheorem{teo}[subsection]{Theorem}

\newtheorem{cor}[subsection]{Corollary}
\newtheorem{prop}[subsection]{Proposition}
\newtheorem{rem}[subsection]{Remark}

\newtheorem{ex}[subsection]{Example}

\let\epsilon\varepsilon
\let\phi=\varphi
\let\kappa=\varkappa
\def\qed{\ifhmode\textqed\fi
   \ifmmode\ifinner\quad\qedsymbol\else\dispqed\fi\fi}
\def\textqed{\unskip\nobreak\penalty50
    \hskip2em\hbox{}\nobreak\hfil\qedsymbol
    \parfillskip=0pt \finalhyphendemerits=0}
\def\dispqed{\rlap{\qquad\qedsymbol}}

\opn\dis{dis}
\def\pnt{{\raise0.5mm\hbox{\large\bf.}}}

\begin{document} 

\title{Generic initial ideal for complete intersections of embedding dimension three with strong Lefschetz property}
\author{Mircea Cimpoea\c{s}}
\address{Mircea Cimpoea\c{s} , Institute of Mathematics of Romanian Academy}
\email{mircea.cimpoeas@imar.ro}

\maketitle
\footnotetext[1]{This paper was supported by the CEEX Program of the Romanian
Ministry of Education and Research, Contract CEX05-D11-11/2005 and by the Higher Education Commission of Pakistan.}

\begin{abstract}
We compute the generic initial ideal of a complete intersection of embedding dimension three with strong Lefschetz
property and we show that it is an almost reverse lexicographic ideal. This enable us to give a proof for Moreno's
conjecture in the case $n=3$.

\vspace{5 pt}
\noindent
\textbf{Keywords:} complete intersection, generic initial ideal, Lefschetz property, Moreno's
conjecture.

\vspace{5 pt}
\noindent
\textbf{2000 Mathematics Subject Classification:} Primary 13P10, Secondary 13D40,13C40.
\end{abstract}

\begin{center}
\textbf{Introduction.}
\end{center}
Let $S=K[x_{1},x_{2},x_{3}]$ be the polynomial ring over a field $K$ of characteristic zero. Let $f_{1}$, $f_{2}$, $f_{3}$ be a regular sequence of homogeneous polynomials of degrees $d_{1},d_{2}$ and $d_{3}$ respectively. We consider the ideal $I=(f_{1},f_{2},f_{3})\subset S$. Obviously, $S/I$ is a complete intersection artinian $K$-algebra. One can easily check that the Hilbert series of $S/I$ depends only on the numbers $d_{1},d_{2}$ and $d_{3}$. More precisely,
\[ H(S/I,t) = (1+t+\cdots +t^{d_{1}-1})(1+t+\cdots +t^{d_{2}-1})(1+t+\cdots +t^{d_{3}-1}).\]
\cite[Lemma 2.9]{HPV} gives an explicit form of $H(S/I,t)$.

We say that a homogeneous polynomial $f$ of degree $d$ is \emph{semiregular} for $S/I$ if the maps 
$(S/I)_{t}\stackrel{\cdot f}{\longrightarrow} (S/I)_{t+d}$ are either injective, either surjective for all 
$t\geq 0$. We say that $S/I$ has the \emph{weak Lefschetz property} (WLP) if there exists a linear form $\ell\in S$,
semiregular on $S/I$. We say that $\ell$ is a weak Lefschetz element for $S/I$. A theorem of Harima-Migliore-Nagel-Watanabe (see \cite{HMNW}) states that $S/I$ has (WLP). We say that $S/I$ has the \emph{strong Lefschetz property} $(SLP)$ if there exists a linear form $\ell\in S$ such that $\ell^{b}$ is semiregular on $S/I$
for all integer $b\geq 1$. In this case, we say that $\ell$ is a strong Lefschetz element for $S/I$. 
Of course, $(SLP)\Rightarrow (WLP)$ but it is not known if $S/I$ has (SLP) for any regular sequence of homogeneous polynomials $f_{1}$, $f_{2}$, $f_{3}$. This is known only in certain cases, for example, when $f_{1},f_{2},f_{3}$ is generic, see \cite{Par} or when $f_{2}\in K[x_{2},x_{3}]$ and $f_{3}\in K[x_{3}]$, see \cite{HW} and \cite{HP}. 

We say that a property $(P)$ holds for a \emph{generic} sequence of homogeneous polynomials 
$f_1,f_2,\ldots, f_n \in S=K[x_1,x_2,\ldots,x_n]$ of given degree $d_1,d_2,\ldots,d_n$ if there exists a nonempty open Zariski subset $U\subset S_{d_1}\times S_{d_2}\times \cdots \times S_{d_n}$ such that for every $(f_{1},f_{2},\ldots,f_{n})\in U$ the property (P) holds. For example, a generic sequence of homogeneous polynomials $f_1,f_2,\ldots,f_n \in S$ is regular.

Now, we present some conjectures and the relations between them (see \cite{Par}). \vspace{5 pt}
\vspace{5 pt}

\noindent
\textbf{Conjecture A}.(Fr\"oberg) If $f_1,f_2,\ldots,f_r\in S=K[x_1,\ldots,x_n]$ is a generic sequence of homogeneous polynomials of given degrees $d_1,d_2,\ldots,d_r$ and $I=(f_1,f_2,\ldots,f_r)$ then the Hilbert series of $S/I$ is
\[ H(S/I) = \left| \frac{\prod_{i=1}^{r} (1-t^{d_i})}{(1-t)^{n}} \right|,\]
where $|\sum_{j\geq 0}a_t t^{j}| = \sum_{j\geq 0}b_t t^{j}$, with $b_{j}=a_{j}$ if $a_{i}>0$ for all $i\leq j$ and $b_{j}=0$ otherwise. \vspace{5 pt}

\noindent
\textbf{Conjecture B}. If $f_1,f_2,\ldots,f_n\in S=K[x_1,\ldots,x_n]$ is a generic sequence of homogeneous polynomials of given degrees $d_1,d_2,\ldots,d_n$ and $I=(f_1,\ldots,f_n)$ then $x_n,x_{n-1},\ldots,x_1$ is a semi-regular sequence on $A=S/I$, i.e. $x_i$ is semiregular on $A/(x_n,\ldots,x_{i+1})$ for all $1\leq i\leq n$. \vspace{5 pt}

\noindent
\textbf{Conjecture C}. If $f_1,f_2,\ldots,f_n\in S=K[x_1,\ldots,x_n]$ is a generic sequence of homogeneous polynomials of given degrees $d_1,d_2,\ldots,d_n$, $I=(f_1,\ldots,f_n)$ and $J$ is the initial ideal of $I$ with respect to the revlex order, then $x_n,x_{n-1},\ldots,x_1$ is a semi-regular sequence on $A=S/(f_1,\ldots,f_n)$. \vspace{5 pt}

\noindent
\textbf{Conjecture D}.(Moreno) If $f_1,f_2,\ldots,f_n\in S=K[x_1,\ldots,x_n]$ is a generic sequence of homogeneous polynomials of given degrees $d_1,d_2,\ldots,d_n$, $I=(f_1,\ldots,f_n)$ and $J$ is the initial ideal of $I$ with respect to the revlex order, then $J$ is an almost revlex ideal, i.e. if $u\in J$ is a minimal generator of $J$ then every monomial of the same degree which preceeds $u$ must be in $J$ as well. \vspace{5 pt}
\vspace{5 pt}

Pardue proved in \cite{Par} that if conjecture $A$ is true for some positive integer $n$ then the conjecture $B$ is true for the same $n$. Also, conjecture $C$ is true for $n$ if and only if $B$ is true for $n$ and if conjecture $B$ is true for some $r$ then $A$ is true for $n<r$ and exactly for that $r$. Also, if conjecture $D$ is true for some $n$ then $B$, and thus $C$, are true for the same $n$. Fr\"oberg \cite{F} and Anick \cite{A} proved that $A$ is true for $n\leq 3$ and so $B$ and $C$ are true for $n\leq 3$. Moreno \cite{Mor} remarked that $D$ is true for $n=2$.

Let $J=Gin_{\leq}(I)$ be the generic initial ideal of $I$, with respect to the reverse lexicographic order. Our aim is to compute $J$ for all regular sequences $f_1,f_2,f_3$ of homogeneous polynomials of given degree $d_1,d_2,d_3$ such 
that $S/I$ has (SLP). We will do this in the sections $2$ and $3$. These computations shown us in particular, that $J$ depends only on the numbers $d_1,d_2,d_3$ (this has been already proved by Popescu and Vladoiu in \cite{HPV}) and more important, that $J$ is an almost reverse lexicographic ideal (Theorem $1.1$). As a consequence, conjecture Moreno ($D$) is true for $n=3$ and $char(K)=0$ (Theorem $1.2$). Note that Conjecture $A$ for $n=3$ and the fact that $J$ is strongly stable does not implies $J$ is almost revlex, and thus does not implies the Moreno's conjecture ($D$) as the example $(x_1^2,x_1x_2,x_2^3,x_1x_3^2)$, shows.

The author wish to thank his Ph.D.adviser, Professor Dorin
Popescu, for support, encouragement and observations on
the content of this paper. Also, he owes a special thank to Dr.Marius Vladoiu for 
his help and for valuables discussions on the subject of this paper.

\section{\textbf{Main results.}}

\begin{teo}
If $f_1,f_2,f_3$ is a regular sequence of homogeneous polynomials of given degrees $d_1,d_2,d_3$ and $I=(f_1,f_2,f_3)$ such that $S/I$ has the (SLP) then $J=Gin(I)$ is uniquely determined and
is an almost reverse lexicographic ideal.
\end{teo}

\begin{proof}
The theorem is a direct consequence of the Propositions $2.3$, $2.8$, $3.3$, $3.8$, $3.13$ and $3.17$.
\end{proof}

\begin{teo}
The conjecture Moreno ($D$) is true for $n=3$ (and $char(K)=0$).
\end{teo}

\begin{proof}
Notice that (SLP) is an open condition. Also, the condition that a sequence of homogeneous polynomial is regular 
is an open condition. It follows, using Theorem $1.1$, that for a generic sequence $f_1,f_2,f_3$ of homogeneous polynomials of given degrees $d_1,d_2,d_3$, $J=Gin(I)$ is almost revlex, where $I=(f_1,f_2,f_3)$. But the definition of the generic initial ideal implies to choose a generic change of variables, and therefore for a generic sequence $f_1,f_2,f_3$ of homogeneous polynomials of given degrees $in(I)$ is almost revlex, as required.
\end{proof}

\begin{rem}
In order to compute the generic initial ideal $J$ we will use the fact that $J$ is a strongly stable ideal, i.e. for any monomial $u\in J$ and any indices $j<i$, if $x_{i}|u$ then $x_{j}u/x_{i}\in J$. Also, a theorem of Wiebe (see \cite{W}) states that $S/I$ has (SLP) if and only if $x_{3}$ is a strong Lefschetz element for $S/J$.
We need to consider several cases: I. $d_1+d_2\leq d_3+1$ with $2$ subcases $d_{1}=d_{2}<d_{3}$, $d_{1}<d_{2}<d_{3}$ (section $2$) and II. $d_1+d_2>d_3+1$ with $4$ subcases: $d_{1}=d_{2}=d_{3}$, $d_{1}=d_{2}<d_{3}$, $d_{1}<d_{2}=d_{3}$, $d_{1}<d_{2}<d_{3}$ (section $3$). 

The construction of $J$ in all cases, follows the next procedure.
For any nonnegative integer $k$, we denote by $J_{k}$ the set of monomials of degree $k$ in $J$. We can easily compute the cardinality of each $J_{k}$ from the Hilbert series of $S/J$. We denote $Shad(J_{k})=\{x_{i}u:\;1\leq i\leq n,\;u\in J_{k}\}$.  We begin with $J_{0}=\emptyset$ and we pass from $J_{k}$ to $J_{k+1}$ noticing that $J_{k+1} = Shad(J_{k})\cup $ eventually some new monomial(s) (exactly $|J_{k+1}|-|Shad(J_{k})|$ new monomials). The fact that $J$ is strongly stable and that $x_{3}$ is a strong Lefschetz element for $S/J$ tell us what we need to add to $Shad(J_{k})$ in order to obtain $J_{k+1}$. We continue this procedure until $k=d_1 + d_2 +d_3 -2$ since $J_{k} = S_{k}$(=the set of all monomials of degree $k$) for any $k\geq d_1 + d_2 +d_3 -2$ and so we cannot add any new monomials in larger degrees. $J$ is the ideal generated by all monomials added to $Shad(J_{k})$ at some step $k$. \vspace{2 pt}
We will present detailed this construction only in the subcase $d_{1}=d_{2}<d_{3}$ of the case $I.d_1+d_2\leq d_3+1$ (Proposition $2.3$), the other cases being presented in sketch, but the reader can easily complete the proofs.
\end{rem}

\newpage
\section{\textbf{Case $d_1+d_2\leq d_3+1$.}}

\medskip
\begin{itemize}
\item Subcase $d_1=d_2<d_3$.
\end{itemize} 
\medskip

\begin{prop}
\em Let $2\leq d:=d_1=d_2<d_3$ be positive integers such that $2d\leq d_3+1$. The Hilbert function of the standard graded complete intersection $A=K[x_{1},x_{2},x_{3}]/I$, where $I$ is the ideal generated by $f_{1}$, $f_{2}$, $f_{3}$, with $f_i$ homogeneous polynomials of degree $d_i$, for all $i$, with $1\leq i\leq 3$, has the form:  
\begin{enumerate}
	\item $H(A,k) = \binom{k+2}{k}$, for $k\leq d-1$.
	\item $H(A,k) = \binom{d+1}{2} +  \sum_{i=1}^{j}(d-i)$, 
	      for $k=d-1+j$, where $0\leq j \leq d-1$.
	\item $H(A,k) = d^{2}$, for $2d-2\leq k \leq d_{3}-1 $.
	\item $H(A,k) = H(A,2d+d_{3}-3-k)$ for $k\geq d_{3}$.
\end{enumerate}
\end{prop}

\begin{proof}
It follows from \cite[Lemma 2.9(a)]{HPV}.
\end{proof}

\begin{cor}
\em In the conditions of Proposition $2.1$, let $J=\Gin(I)$ be the generic initial ideal of $I$ with respect to the reverse lexicographic order. If we denote by $J_{k}$ the set of monomials of $J$ of degree $k$, then:
\begin{enumerate}
	\item $J_{k} = \emptyset$, for $k\leq d-1$.
	\item $|J_{k}| = j(j+1)$, for $k=d-1+j$, where $0\leq j \leq d-1$.
	\item $|J_{k}| = d(d-1) + 2dj + \frac{j(j-1)}{2}$, for $k=2d-2+j$, where $0\leq j \leq d_3+1-2d$.
  \item $|J_{k}| = \frac{d_{3}(d_{3}+1)}{2} -d^2 +jd_{3} + j(j+1)$, for $k=d_{3}-1+j$, where $0\leq j \leq d-1$.
  \item $|J_{k}| = \frac{d_3(d_3-1)}{2}+d(d_3-1)+j(d_{3}+2d)$, for $k=d+d_{3}-2+j$, where $0\leq j \leq d$.
  \item $J_{k} = S_{k}$, for $k\geq 2d+d_{3}-2$, where $S_{k}$ is the set of monomials of degree $k$.
\end{enumerate}
\end{cor}

\begin{proof}
Using that $|J_k| = |S_k| - H(S/J,k)$, together with the general fact that $H(S/J,k)=H(S/I,k)$, the proof follows immediately from Proposition $2.1$. 
\end{proof}

\begin{prop}
\em Let $2\leq d:=d_1=d_2<d_3$ be positive integers such that $2d\leq d_3+1$. Let $f_{1},f_{2},f_{3}\in K[x_{1},x_{2},x_{3}]$ be a regular sequence of homogeneous polynomials of degrees $d_{1},d_{2},d_{3}$. If  $I=(f_{1},f_{2},f_{3})$, and $J=\Gin(I)$, the generic initial ideal with respect to the reverse lexicographic order, and
$S/I$ has (SLP), then:
\[ J = ( x_{1}^{d}, x_{1}^{d-j-1}x_{2}^{2j+1}\;for\;0\leq j\leq d-1, 
x_{2}^{2d-2j-2}x_{3}^{d_{3}-2d+2j+2}\{x_{1},x_{2}\}^{j}\;for\;0\leq j\leq d-2,\]
\[ x_{3}^{d_{3}+2j-2}\{x_{1},x_{2}\}^{d-j}\; for\; 1\leq j\leq d). \]
\end{prop}

\begin{proof}
We have $|J_{d}|=2$, hence 
$J_{d}=\{ x_{1}^{d-1}\{x_{1},x_{2}\} \}$, since $J$ is a strongly stable ideal. Therefore:
\[ 
Shad(J_{d}) = \{x_{1}^{d-1}\{x_{1},x_{2}\}^{2} ,x_{1}^{d-1}x_{3}\{x_{1},x_{2}\}\},
\]
 Now we have two possibilities to analyze: $d=2$ and $d\geq 3$. First, suppose $d\geq 3$.

Using the formulae from Corollary $2.2$ we have $|J_{d+1}|-|Shad(J_{d})|= 1$, so there is only one generator to add to the set $Shad(J_{d})$ in order to obtain $J_{d+1}$.  Since $J$ is strongly stable, we have only two possibilities: $x_{1}^{d-2}x_{2}^{3}$ or $x_{1}^{d-1}x_{3}^{2}$. We cannot have $x_{1}^{d-1}x_{3}^{2}$, since otherwise the application
\[ 
(S/J)_{d-1} \stackrel{\cdot x_{3}^{2}}{\longrightarrow} (S/J)_{d+1},
\]
with $|(S/J)_{d-1}|<|(S/J)_{d+1}|$ (see Proposition $2.1$) would not be injective ($0\neq x_{1}^{d-1}\in (S/J)_{d-1}$ and is mapped to $0$), which is a contradiction to the fact that $x_{3}$ is a strong Lefschetz element for $S/J$. Hence:
\[
J_{d+1} = \{x_{1}^{d-2}\{x_{1},x_{2}\}^{3} ,x_{1}^{d-1}x_{3}\{x_{1},x_{2}\}\}.
\]
We prove by induction on $j$, with $1\leq j \leq d-2$, that:
\[ 
J_{d+j} = Shad(J_{d+j-1})\cup \{x_{1}^{d-j-1}x_{2}^{2j+1}\} =
\]
\[
= \{x_{1}^{d-j-1}\{x_{1},x_{2}\}^{2j+1}, x_{1}^{d-j}x_{3}\{x_{1},x_{2}\}^{2j-1},\ldots,x_{1}^{d-1}x_{3}^{j}\{x_{1},x_{2}\}\}.
\]
The assertion was checked above for $j=1$. Assume now that the statement is true for some $j<d-2$. Then $Shad(J_{d+j})$ is the following set:
\[
\{x_{1}^{d-j-1}\{x_{1},x_{2}\}^{2j+2},x_{1}^{d-j-1}x_{3}\{x_{1},x_{2}\}^{2j+1}, x_{1}^{d-j}x_{3}^2\{x_{1},x_{2}\}^{2j-1},\ldots,x_{1}^{d-1}x_{3}^{j+1}\{x_{1},x_{2}\}\}. 
\]
We have $|J_{d+j+1}| - |Shad(J_{d+j})| = 1$, so we must add only one generator to $Shad(J_{d+j})$ to get $J_{d+j+1}$. The ideal $J$, being strongly stable, allows only two possibilities, namely $x_{1}^{d-j-2}x_{2}^{2j+3}$ or $x_{1}^{d-2}x_{2}^2x_3^{j+1}$. The second one is not allowed because the application
\[
(S/J)_{d} \stackrel{\cdot x_{3}^{j+1}}{\longrightarrow} (S/J)_{d+j+1}
\]
would not be injective ($0\neq x_{1}^{d-2}x_2^2\in (S/J)_{d}$ and is mapped to $0$), $x_{3}$ being a strong-Lefschetz element for $S/J$. Therefore, we must add $x_{1}^{d-j-2}x_{2}^{2j+3}$ and our claim is proved. In particular, we obtain
\[
J_{2d-2} =\{x_{1}\{x_{1},x_{2}\}^{2d-3}, x_{1}^{2}x_{3}\{x_{1},x_{2}\}^{2d-5},\ldots,x_{1}^{d-1}x_{3}^{d-2}\{x_{1},x_{2}\}\}.
\]

In order to compute $J_{k}$, with $2d-2\leq k \leq d_3$, we must consider two posibilities. 

\begin{itemize}
	\item 1. $d_3=2d-1$. 
\end{itemize}

Since $|J_{2d-1}|-|Shad(J_{2d-2})| = 2$, there are two generators to add to $Shad(J_{2d-2})$. We prove that these generators are
$x_2^{2d-1},x_2^{2d-2}x_3$. Assuming by contradiction that we have other generators, since $J$ is strongly stable, it follows that
there is at least one generator from the set $\{x_{1}x_{2}^{2d-4}x_{3}^{2},x_{1}^{2}x_{2}^{2d-6}x_{3}^{3},\ldots, x_{1}^{d-2}x_{2}^{2}x_{3}^{d-1}\}$. Then, the application $(S/J)_{2d-3} \stackrel{\cdot x_{3}^{2}}{\longrightarrow} (S/J)_{2d-1}$
would not be injective, a contradiction since $x_3$ is a strong Lefschetz element for $S/J$. Therefore
\[
J_{2d-1} = J_{d_3} = \{ \{x_{1},x_{2}\}^{2d-1},  x_{3}\{x_{1},x_{2}\}^{2d - 2},\ldots,x_{1}^{d-1}x_{3}^{d-1}\{x_{1},x_{2}\}\}.
\]

\begin{itemize}
	\item 2. $d_3>2d-1$. 
\end{itemize}	

Since $|J_{2d-1}|-|Shad(J_{2d-2})| = 1$, there is only one generator to add to $Shad(J_{2d-2})$, which can be selected from the set $\{x_2^{2d-1},x_1x_2^{2d-4}x_3^2,x_1^2x_2^{2d-6}x_3^3,\ldots,x_1^{d-2}x_2^2x_3^{d-1}\}$ because $J$ is strongly stable. In a similar manner to what we have done above can be shown that, $x_3$ being a strong Lefschetz  element, leaves us as unique possibility $x_{2}^{2d-1}$, therefore:
\[
   J_{2d-1}=\{\{x_{1},x_{2}\}^{2d-1}, x_{1}x_{3}\{x_{1},x_{2}\}^{2d-3},\ldots,x_{1}^{d-1}x_{3}^{d-1}\{x_{1},x_{2}\}\}. 
\] 
One can easily show, using induction on $1\leq j \leq d_{3}-2d$, if case, that $|J_{2d-1+j}| = |Shad(J_{2d-2+j})|$ and $J_{2d-1+j}$ is the set \[ \{ \{x_{1},x_{2}\}^{2d-1+j}, \ldots ,x_{3}^{j}\{x_{1},x_{2}\}^{2d-1},
x_{3}^{j+1}x_{1}\{x_{1},x_{2}\}^{2d-2},\ldots,x_{3}^{d+j-1}x_{1}^{d-1}\{x_{1},x_{2}\}\} .\]
In particular, we obtain that $J_{d_{3}-1}$ is the set
\[ \{\{x_{1},x_{2}\}^{d_{3}-1}, \ldots,x_{3}^{d_{3}-2d}\{x_{1},x_{2}\}^{2d-1}, x_{3}^{d_{3}-2d+1}x_{1}\{x_{1},x_{2}\}^{2d-3},\ldots,x_{3}^{d_{3}-d-1}x_{1}^{d-1}\{x_{1},x_{2}\}\}.\]

Since $|J_{d_3}|-|Shad(J_{d_{3}-1})|=1$, the generator which has to be add to $Shad(J_{d_{3}-1})$ can be selected from the set
$\{x_{2}^{2d-2}x_3^{d_3 - 2d+2}, x_{1}x_2^{2d-4}x_3^{d_3 - 2d+3}, \ldots ,x_{1}^{d-2}x_{2}^{2}x_{3}^{d_{3}-d}\}$ such that $J$
is strongly stable. The generator is $x_{2}^{2d-2}x_3^{d_3 - 2d+2}$, otherwise the application $(S/J)_{2d-3} \stackrel{\cdot x_{3}^{d_{3}-2d+3}}{\longrightarrow} (S/J)_{d_{3}}$ is not injective, a contradiction, since $x_3$ is a strong Lefschetz element
for $S/J$. Hence, we get that $J_{d_{3}}$ is
\[  \{\{x_{1},x_{2}\}^{d_{3}}, \ldots, x_{3}^{d_{3}-2d+2}\{x_{1},x_{2}\}^{2d-2},x_{1}^{2}x_{3}^{d_{3}-2d+3}\{x_{1},x_{2}\}^{2d-5},\ldots,x_{3}^{d_{3}-d}x_{1}^{d-1}\{x_{1},x_{2}\}\},\] and one can check that is the same formula as in 1.($d_3 = 2d-1$).


Now, we show by induction on $1\leq j\leq d-2$ that
\[ 
J_{d_{3}+j} = Shad(J_{d_{3}-1+j}) \cup \{x_{2}^{2d-2j-2}x_{3}^{d_{3}-2d+2j+2}\{x_{1},x_{2}\}^{j} \}. 
\]
Indeed, for $j=1$, $|J_{d_{3}+1}|-|Shad(J_{d_{3}})|=2$ and the generators which must be added are $x_{1}x_{2}^{2d-4}x_{3}^{d_{3}-2d+4}, x_{2}^{2d-3}x_{3}^{d_{3}-2d+4}$. If not, since $J$ is strongly stable, then at least one of the generators belongs to the set
$\{ x_{1}^{2}x_{2}^{2d-6}x_{3}^{d_{3}-2d+5},\ldots,x_{1}^{d-2}x_{2}^{2}x_{3}^{d_{3}-d+1} \}$ (for $d=3$ this is the emptyset).
but then the map $(S/J)_{2d-4} \stackrel{x_{3}^{d_{3}-2d+5}}{\longrightarrow} (S/J)_{d_{3}+1}$ is not injective, contradiction.

Assume now that we proved the assertion for some $j<d-2$. Then $|J_{d_{3}+j+1}|-|Shad(J_{d_{3}+j})| = j+2$ and the new generators are $x_{2}^{2d-2j-4}x_{3}^{d_{3}-2d+2j+4}\{x_{1},x_{2}\}^{j+1}$. Indeed, if not, since $J$ is strongly stable, then at least one of the generators belongs to the set
$\{x_{1}^{j+2}x_{2}^{2d-2j-6}x_{3}^{d_{3}-2d+2j+5},\ldots,x_{1}^{d-2}x_{2}^{2}x_{3}^{ d_{3}-d+j+1}\}$ (for $d=3$ this is the emptyset...) but then the map $(S/J)_{2d-j-4} \stackrel{x_{3}^{d_{3}-2d+2j+5}}{\longrightarrow} (S/J)_{d_{3}+j+1}$ is not injective, contradiction, and we are done. Hence, 
\[J_{d_{3}+d-2} = \{ \{x_{1},x_{2}\}^{d_{3}+d-2},x_{3}\{x_{1},x_{2}\}^{d+d_{3}-3},\ldots ,x_{3}^{d_{3}-2}\{x_{1},x_{2}\}^{d}\}.\]

We prove by induction on $1\leq j \leq d$ that $J_{d+d_{3}-2+j} = Shad(J_{d+d_{3}-3+j})\cup x_{3}^{d_{3}+2j-2}\{x_{1},x_{2}\}^{d-j}$.
If $j=1$ then $|J_{d+d_{3}-1}|- |Shad(J_{d+d_{3}-2})| = d$ so we must add $d$ generators, which are precisely the elements of the set $x_{3}^{d_{3}}\{x_{1},x_{2}\}^{d-1}$. Indeed, if we have a generator which does not belong to the set it is divisible by $x_{3}^{d_{3}+1}$ and therefore the map $(S/J)_{d-2} \stackrel{\cdot x_{3}^{d_{3}+1}}{\longrightarrow} (S/J)_{d_{3}+d-1}$ is not injective, which is a contradiction with $x_{3}$ is a strong Lefschetz element for $S/J$ (the map has to be bijective). The induction
step is similar and finally we obtain that $J_{d_{3}+2d-2}=S_{d_{3}+2d-2}$ and thus we cannot add new minimal generators of $J$ in degree $>d_{3}+2d-2$.

In order to complete the proof we must consider now $d=2$. The hypothesis implies $d_{3}\geq 3$. We already seen that $J_{2}=\{x_{1}^{2},x_{1}x_{2}\}$ and $Shad(J_{2}) = \{x_{1}\{x_{1},x_{2}\}^{2} ,x_{1}x_{3}\{x_{1},x_{2}\}\}$.

Using the formulae from Corollary $2.2$ we have $|J_{3}|-|Shad(J_{2})|= 1$, so there is only one generator to add to the set $Shad(J_{d})$ in order to obtain $J_{d+1}$. 

Since $J$ is strongly stable, we have only two possibilities: $x_{2}^{3}$ or $x_{1}x_{3}^{2}$. We can not have $x_{1}x_{3}^{2}$, since otherwise the application
\[ 
(S/J)_{1} \stackrel{\cdot x_{3}^{2}}{\longrightarrow} (S/J)_{3},
\]
with $|(S/J)_{1}|<|(S/J)_{3}|$ (see Proposition $2.1$) would not be injective ($0\neq x_{1}\in (S/J)_{1}$ and is mapped to $0$), which is a contradiction to the fact that $x_{3}$ is a strong Lefschetz element for $S/J$. Hence:
\[
J_{3} = \{\{x_{1},x_{2}\}^{3} ,x_{1}x_{3}\{x_{1},x_{2}\}\}.
\]
Assume now $d_{3}\geq 4$. One can easily show, using induction on $1\leq j \leq d_{3}-4$, if case, that $|J_{3+j}| = |Shad(J_{2+j})|$ and $J_{3+j}$ is the set \[ \{ \{x_{1},x_{2}\}^{3+j}, \ldots ,x_{3}^{j}\{x_{1},x_{2}\}^{3},
x_{3}^{j+1}x_{1}\{x_{1},x_{2}\}\}.\]
In particular, we obtain that \[J_{d_{3}-1} = \{\{x_{1},x_{2}\}^{d_{3}-1}, \ldots,x_{3}^{d_{3}-4}\{x_{1},x_{2}\}^{3}, x_{3}^{d_{3}-3}x_{1}\{x_{1},x_{2}\}\}.\]

Since $|J_{d_3}|-|Shad(J_{d_{3}-1})|=1$, the generator which has to be add to $Shad(J_{d_{3}-1})$ is exactly $x_{2}^{2}x_3^{d_3 - 2}$ such that $J$ is strongly stable. Hence, we get
\[J_{d_{3}} =  \{\{x_{1},x_{2}\}^{d_{3}}, \ldots, x_{3}^{d_{3}-3}\{x_{1},x_{2}\}^{3} , x_{1}x_{3}^{d_{3}-2}\{x_{1},x_{2}\}\},\] and one can check that is the same formula as in the case $d_3 = 3$.

Since $|J_{d_3+1}|-|Shad(J_{d_{3}})|=1$, there is only one generator to add to the set $Shad(J_{d_{e}})$ in order to obtain $J_{d_3+1}$. Since $J$ is strongly stable, we have only two possibilities: $x_{2}^{2}x_{3}^{d_{3}-1}$ or $x_{1}x_{3}^{d_{3}}$. We can not have $x_{1}x_{3}^{d_{3}}$, since otherwise the application
\[ 
(S/J)_{1} \stackrel{\cdot x_{3}^{d_{3}}}{\longrightarrow} (S/J)_{d_{3}+1},
\]
with $|(S/J)_{1}|=|(S/J)_{3}|$ (see Proposition $2.1$) would not be injective ($0\neq x_{1}\in (S/J)_{1}$ and is mapped to $0$), which is a contradiction to the fact that $x_{3}$ is a strong Lefschetz element for $S/J$. Hence:
\[
J_{d_{3}+1} = \{\{x_{1},x_{2}\}^{d_{3}+1}, \ldots, x_{3}^{d_{3}-1}\{x_{1},x_{2}\}^{2} \}.
\] 
Since $|J_{d_3+2}|-|Shad(J_{d_{3}+1})|=2$ and $J$ is strongly stable, we must add $x_{1}x_{3}^{d_{3}+1}$ and
$x_{2}x_{3}^{d_{3}+1}$ at $Shad(J_{d_{3}+1})$ in order to obtain $J_{d_3+2}$. Hence 
$J_{d_{3}+2}=S_{d_{3}+2}\setminus \{x_{3}^{d_{3}+2}\}$. Finally, since $J_{d_{3}+3}=S_{d_{3}+3}$ we add $x_{3}^{d_{3}+3}$ at $Shad(J_{d_{3}+2})$ and thus we cannot add new minimal generators of $J$ in degree $>d_{3}+2$.
\end{proof}

\begin{cor}
In the conditions of the above proposition, the number of minimal generators of $J$ is $d^{2}+d+1$.
\end{cor}

\begin{ex}
Let $d_{1}=d_{2}=3$ and $d_{3}=9$. Proposition $2.3$ implies:
\[ J=(x_{1}^{3},\; x_{1}^{2}x_{2},\; x_{1}x_{2}^{3},\; x_{2}^{5},\; x_{2}^{4}x_{3}^{5},\; x_{1}x_{2}^{2}x_{3}^{7},\;
      x_{2}^{3}x_{3}^{7},\;x_{3}^{9}\{x_{1},x_{2}\}^{2},\;x_{3}^{11}\{x_{1},x_{2}\},\;x_{3}^{13}    ). \]
\end{ex}

\newpage
\medskip
\begin{itemize}
\item Subcase $d_1<d_2<d_3$.
\end{itemize} 
\medskip

\begin{prop}
\em Let $2\leq d_1<d_2<d_3$ be positive integers such that \linebreak $d_1 + d_2\leq d_3+1$. The Hilbert function of the standard graded complete intersection $A=K[x_{1},x_{2},x_{3}]/I$, where $I$ is the ideal generated by $f_{1}$, $f_{2}$, $f_{3}$, with $f_i$ homogeneous polynomials of degree $d_i$, for all $i$, with $1\leq i\leq 3$, has the form:  
\begin{enumerate}
	\item $H(A,k) = \binom{k+2}{k}$, for $k\leq d_{1}-1$.
	\item $H(A,k) = \binom{d_{1}+1}{2} + jd_{1}$, for $k=j + d_1 - 1$, where $0\leq j\leq d_2 - d_1$.
	\item $H(A,k) = \binom{d_{1}+1}{2} + d_{1}(d_{2}-d_{1}) + \sum_{i=1}^{j}(d_{1}-i)$, 	      
	      for $k=j+d_{2}-1$, where $0 \leq j \leq d_{1}-1$.
	\item $H(A,k) = d_{1}d_{2}$, for $d_{1}+d_{2}-2\leq k \leq d_{3}-1 $.
	\item $H(A,k) = H(A,d_{1}+d_{2}+d_{3}-3-k)$ for $k\geq d_{3}$.
\end{enumerate}
\end{prop}

\begin{proof}
It follows from \cite[Lemma 2.9(a)]{HPV}.
\end{proof}

\begin{cor}
\em In the conditions of Proposition $2.6$, let $J=\Gin(I)$ be the generic initial ideal of $I$ with respect to the reverse lexicographic order. If we denote by $J_{k}$ the set of monomials of $J$ of degree $k$, then:
\begin{enumerate}
	\item $|J_{k}| = 0$, for $k\leq d-1$.
	\item $|J_{k}| = j(j+1)/2$, for $k=j+d_1-1$, where $0\leq j\leq d_2 - d_1$.	
	\item $|J_{k}| = \frac{(d_{2}-d_{1})((d_{2}-d_{1}-1))}{2} + j(d_{2}-d_{1}) + j(j+1)$, 
	      for $k=j+d_{2}-1$, where $0 \leq j \leq d_{1}-1$.	      
  \item $|J_{k}| = \frac{d_{1}^{2}+d_{2}^{2}-d_{1}-d_{2}}{2} + j(d_{1}+d_{2}) + \frac{j(j-1)}{2}$, 
        for $k=j+d_{1}+d_{2}-2$, where $0 \leq j \leq d_{3}-d_{1}-d_{2}+1$.        
  \item $|J_{k}| = \frac{d_{3}^{2}+d_{3}-2d_{1}d_{2}}{2} +jd_{3} + j(j+1)$, 
        for $k=j+d_{3}-1$ where $0 \leq j \leq d_{1}-1$.
  
  \item $|J_{k}| = \frac{ (d_{1}+d_{3})(d_{1}+d_{3}-1) +d_{1}^{2} - d_{1} -2d_{1}d_{2}}{2} +j(d_{3}+2d_{1}) +
         \frac{j(j-1)}{2}$, for $k = j + d_1 + d_3 - 2$, where $0 \leq j \leq d_{2}-d_{1}$ .
  \item $|J_{k}| = \frac{(d_{2}+d_{3})(d_{2}+d_{3}-1) +d_{1}(d_{1}-1)}{2} + j(d_{1}+d_{2}+d_{3})$, 
        for $k=d_{2}+d_{3}-2$, where $0 \leq j \leq d_{1}-1$.	
  \item $J_{k} = S_{k}$, for $k\geq 3d-2$.
\end{enumerate}
\end{cor}

\begin{prop}
\em Let $2\leq d_1<d_2<d_3$ be positive integers such that $d_1 + d_2\leq d_3+1$. Let $f_{1},f_{2},f_{3}\in K[x_{1},x_{2},x_{3}]$ be a regular sequence of homogeneous polynomials of degrees $d_{1},d_{2},d_{3}$. If $I=(f_{1},f_{2},f_{3})$ ,$J=\Gin(I)$, the generic initial ideal with respect to the reverse lexicographic order, and $S/I$ has (SLP),then:
\[ J = ( x_{1}^{d_{1}}, x_{1}^{d_{1}-j} x_{2}^{d_{2}-d_{1}+2j-1}\;for\;1\leq j\leq d_1 - 1,\; x_{2}^{d_1+d_2 -1},\;
x_{2}^{d_1+d_2 -2}x_{3}^{d_{3}-d_1-d_2+2}, \]
\[ x_{3}^{d_{3}-d_{1}-d_{2}+2j+2}x_{2}^{d_{1}+d_{2}-2j-2}\{x_{1},x_{2}\}^{j}\;for\;1\leq j\leq d_1-2,\; \]\[
x_{3}^{d_{3}+d_{1}-d_{2}-2+2j}x_{2}^{d_{2}-d_{1}+1-j}\{x_{1},x_{2}\}^{d_1 - 1}\;for\;1\leq j\leq d_2-d_1,\;\]\[
x_{3}^{d_{3}+d_2-d_1+2j-2}\{x_{1},x_{2}\}^{d_1-j}\;for\;1\leq j\leq d_1 ).\]
\end{prop}

\begin{proof}
We have $|J_{d_1}|=1$, hence $J_{d_1}=\{x_{1}^{d_1}\}$, since $J$ is a strongly stable ideal. Therefore:
\[ 
Shad(J_{d_1}) = \{x_{1}^{d_1}\{x_{1},x_{2}\} ,x_{1}^{d_1}x_{3}\}.
\]
Assume $d_2>d_1+1$. Since $|J_{d_1+1}|=|Shad(J_{d_1})|$ from the formulae of $2.7$, if follows $J_{d_1+1}=Shad(J_{d_1})$. We prove by induction on $1\leq j\leq d_2-d_1-1$ that 
\[ J_{d_1+j}=Shad(J_{d_1+j-1})= \{x_{1}^{d_{1}}\{x_{1},x_{2}\}^{j}, x_{3}x_{1}^{d_{1}}\{x_{1},x_{2}\}^{j-1},\ldots,x_{3}^{j}x_{1}^{d_{1}} \}. \]
Indeed, the case $j=1$ is already proved. Suppose the assertion is true for some $j<d_2-d_1-1$. Since
$|J_{d_1+j+1}|-|Shad(J_{d_1+j})|=0$ it follows that 
\[ J_{d_1+j+1}=Shad(J_{d_1+j})= \{x_{1}^{d_{1}}\{x_{1},x_{2}\}^{j+1}, x_{3}x_{1}^{d_{1}}\{x_{1},x_{2}\}^{j},\ldots,x_{3}^{j+1}x_{1}^{d_{1}} \}\]
thus we are done. In particular, we get
\[ J_{d_{2}-1} = \{x_{1}^{d_{1}}\{x_{1},x_{2}\}^{d_{2}-d_{1}-1}, x_{3}x_{1}^{d_{1}}\{x_{1},x_{2}\}^{d_{2}-d_{1}-2},\ldots,x_{3}^{d_{2}-d_{1}-1}x_{1}^{d_{1}}\} \]
which is the same formula as in the case $d_2=d_1+1$.

We have $|J_{d_2}| - |Shad(J_{d_{2}-1})| = 1$ so we must add a new generator to $Shad(J_{d_{2}-1})$ to obtain $J_{d_{2}}$. Since $J$ is strongly stable and $x_3$ is a strong Lefschetz element for $S/J$ this new generator is $x_{1}^{d_{1}-1}x_{2}^{d_{2}-d_{1}+1}$, therefore
\[ J_{d_{2}} = \{x_{1}^{d_{1}-1}\{x_{1},x_{2}\}^{d_{2}-d_{1}+1}, x_{3}x_{1}^{d_{1}}\{x_{1},x_{2}\}^{d_{2}-d_{1}-1},\ldots,x_{3}^{d_{2}-d_{1}}x_{1}^{d_{1}}\}.\]
Assume $d_1>2$. We prove by induction on $1\leq j\leq d_{1}-1$ that
\[ J_{d_{2}-1+j} = Shad(J_{d_{2}-2+j}) \cup \{ x_{1}^{d_{1}-j} x_{2}^{d_{2}-d_{1}+2j-1}\} =  \{x_{1}^{d_{1}-j}\{x_{1},x_{2}\}^{d_{2}-d_{1}+2j-1},\]\[ x_{3}x_{1}^{d_{1}-j+1}\{x_{1},x_{2}\}^{d_{2}-d_{1}+2j-3}, \ldots, x_{3}^{j}x_{1}^{d_{1}}\{x_{1},x_{2}\}^{d_{2}-d_{1}-1},\ldots,x_{3}^{d_{2}-d_{1}+j-1}x_{1}^{d_{1}}\}.\]
The assertion was proved for $j=1$. Suppose $1\leq j<d_1-1$ and the assertion is true for $j$. We have
$|J_{d_{2}+j}|-|Shad(J_{d_{2}-1+j})|=1$, thus we must add a new generator to $Shad(J_{d_{2}-1+j})$ in order to
obtain $J_{d_{2}+j}$ and since $J$ is strongly stable and $x_3$ is a strong Lefschetz element for $S/J$, this is $x_{1}^{d_{1}-j-1} x_{2}^{d_{2}-d_{1}+2j+1}$ and we are done. In particular, we obtain:
\[ J_{d_{1}+d_{2}-2} = \{x_{1}\{x_{1},x_{2}\}^{d_{1}+d_{2}-3}, x_{3}x_{1}^{2}\{x_{1},x_{2}\}^{d_{1}+d_{2}-5}, \ldots,\] 
\[ x_{3}^{d_{1}-1}x_{1}^{d_{1}}\{x_{1},x_{2}\}^{d_{2}-d_{1}-1},
x_{3}^{d_{1}}x_{1}^{d_{1}}\{x_{1},x_{2}\}^{d_{2}-d_{1}-2},\ldots,x_{3}^{d_{2}-2}x_{1}^{d_{1}} \}\]
and one can check that is the same expression as in the case $d_1=2$.

In order to compute $J_{k}$, with $2d-2\leq k \leq d_3$, we must consider two possibilities. 

\begin{itemize}
	\item 1.$d_3=d_1+d_2-1$. 
\end{itemize}

Since $|J_{d_1+d_2-1}|-|Shad(J_{d_1+d_2-2})| = 2$, there are two generators to add to $Shad(J_{d_1+d_2-2})$ to get $J_{d_1+d_2-1}$, but on the other hand $J$ is strongly stable and $x_3$ is a strong Lefschetz element for $S/J$ so these generators must be $x_2^{d_1+d_2-1},x_2^{d_1+d_2-2}x_3$. Therefore:
\[J_{d_1+d_2-1} = J_{d_3} = \{ \{x_{1},x_{2}\}^{d_1+d_2-1},  x_{3}\{x_{1},x_{2}\}^{d_1 +d_2 - 2}, \ldots,x_{1}^{d_1}x_{3}^{d_2-1}\}.\] \pagebreak

\begin{itemize}
	\item 2.$d_3>d_1+d_2-1$. 
\end{itemize}	

Since $|J_{d_1+d_2-1}|-|Shad(J_{d_1+d_2-2})| = 1$, there is only one generator to add to $Shad(J_{2d-2})$, which is precisely $x_{2}^{d_1+d_2-1}$ since $J$ is strongly stable and $x_3$ is a strong Lefschetz element for $S/J$. Therefore
\[ J_{d_1+d_2-1} = \{\{x_{1},x_{2}\}^{d_1+d_2-1}, x_{1}x_{3}\{x_{1},x_{2}\}^{d_1+d_2-3},\ldots, x_{1}^{d_1}x_{3}^{d_2-1} \}.\] 
One can easily show, using induction on $1\leq j \leq d_{3}-d_1-d_2$, if case, that \linebreak
$|J_{d_1+d_2-1+j}| = |Shad(J_{d_1+d_2-2+j})|$ and $J_{d_1+d_2-1+j}$ is the set 
\[\{ \{x_{1},x_{2}\}^{j+d_{1}+d_{2}-1}, x_{3}\{x_{1},x_{2}\}^{j+d_{1}+d_{2}-2}, \ldots, x_{3}^{j}\{x_{1},x_{2}\}^{d_{1}+d_{2}-1},\]
\[ x_{3}^{j+1}x_{1}\{x_{1},x_{2}\}^{d_{1}+d_{2}-3},\ldots,x_{3}^{d_{1}+j}x_{1}^{d_{1}}\{x_{1},x_{2}\}^{d_{2}-d_{1}-1},
\ldots, x_{3}^{j+d_{2}-1}x_{1}^{d_{1}} \} \]

\[So\;J_{d_{3}-1}=\{ \{x_{1},x_{2}\}^{d_{3}-1}, x_{3}\{x_{1},x_{2}\}^{d_{3}-2}, \ldots, x_{3}^{d_{3}-d_{1}-d_{2}}\{x_{1},x_{2}\}^{d_{1}+d_{2}-1},\]
\[ x_{3}^{d_{3}-d_{1}-d_{2}+1}x_{1}\{x_{1},x_{2}\}^{d_{1}+d_{2}-3},\ldots,
x_{3}^{d_{3}-d_{2}}x_{1}^{d_{1}}\{x_{1},x_{2}\}^{d_{2}-d_{1}-1},
\ldots, x_{3}^{d_{3}-d_{1}-1}x_{1}^{d_{1}} \} \]

Since $|J_{d_3}|-|Shad(J_{d_{3}-1})|=1$, $J$ is strongly stable and $x_3$ is a strong Lefschetz element for $S/J$, the generator which has to be added to $Shad(J_{d_{3}-1})$ is $x_{2}^{d_1+d_2-2}x_3^{d_3 - d_1-d_2 +2}$. Hence, we get 
\[ J_{d_3} = \{ \{x_{1},x_{2}\}^{d_{3}}, x_{3}\{x_{1},x_{2}\}^{d_{3}-1}, \ldots, x_{3}^{d_{3}-d_{1}-d_{2}+2}\{x_{1},x_{2}\}^{d_{1}+d_{2}-2}, \]\[
x_{3}^{d_{3}-d_{1}-d_{2}+3}x_{1}^{2}\{x_{1},x_{2}\}^{d_{1}+d_{2}-5},\ldots,
x_{3}^{d_{3}-d_{2}+1}x_{1}^{d_{1}}\{x_{1},x_{2}\}^{d_{2}-d_{1}-1},
\ldots, x_{3}^{d_{3}-d_{1}}x_{1}^{d_{1}}, \} \]
and one can check that is the same formula as in 1.($d_3 = d_1+d_2-1$).

Assume now $d_1>2$. 
We show by induction on $1\leq j \leq d_1-2$ that
\[ J_{d_3+j} = Shad(J_{d_3-1+j})\cup \{x_{1}^{j}x_{3}^{d_{3}-d_{1}-d_{2}+2j+2}x_{2}^{d_{1}+d_{2}-2j-2}, \ldots,   x_{3}^{d_{3}-d_{1}-d_{2}+2j+2}x_{2}^{d_{1}+d_{2}-2-j}\}. \]
Indeed, for $j=1$, $|J_{d_3+1}|-|Shad(J_{d_3})|=2$ so we must add two generators to $Shad(J_{d_3})$ in order to
obtain $J_{d_3+1}$. Since $J$ is strongly stable and $x_3$ is a strong Lefschetz element for $S/J$, these new generators are $x_{1}x_{3}^{d_{3}-d_{1}-d_{2}+4}x_{2}^{d_{1}+d_{2}-2j-4}$ and $x_{3}^{d_{3}-d_{1}-d_{2}+4}x_{2}^{d_{1}+d_{2}-2j-3}$.
Assume now that we proved the assertion for some $j<d_1-2$. Then $|J_{d_3+j+1}|-|Shad(J_{d_3+j})| = j+2$ and since
$J$ is strongly stable and $x_3$ is a strong Lefschetz element for $S/J$, the new generators are $x_{3}^{d_{3}-d_{1}-d_{2}+2j+4}x_{2}^{d_{1}+d_{2}-2j-4}\{x_{1},x_{2}\}^{j+1}$ as required. Hence, we get
\[ J_{d_{1}+d_{3}-2} = \{ \{x_{1},x_{2}\}^{d_{1}+d_{3}-2},\ldots,x_{3}^{d_{3}+d_{1}-d_{2}-2}\{x_{1},x_{2}\}^{d_{2}},\]
\[ x_{1}^{d_{1}}x_{3}^{d_{3}+d_{1}-d_{2}-1}\{x_{1},x_{2}\}^{d_{2}},\ldots,x_{1}^{d_{1}}x_{3}^{d_{3}-2}\},\]
and one can check that is the same formula as in the case $d_1 = 2$.

We have $|J_{d_{1}+d_{3}-1}|-|Shad( J_{d_{1}+d_{3}-2})|=d_1$ so we must add $d_1$ new generators to
$Shad( J_{d_{1}+d_{3}-2})$ and since $J$ is strongly stable and $x_3$ is a strong Lefschetz element for $S/J$, they are $x_{3}^{d_{3}+d_{1}-d_{2}}x_{2}^{d_{2}-d_{1}}\{x_{1},x_{2}\}^{d_1 - 1}$. Therefore $J_{d_{1}+d_{3}-1}$ is the set
\[\{ \{x_{1},x_{2}\}^{d_{1}+d_{3}-1}, .., x_{3}^{d_{3}+d_{1}-d_{2}}\{x_{1},x_{2}\}^{d_{2}-1},
x_{1}^{d_{1}}x_{3}^{d_{3}+d_{1}-d_{2}+1}\{x_{1},x_{2}\}^{d_{2}-d_{1}-2},.., x_{1}^{d_{1}}x_{3}^{d_{3}-1}\}. \]

Suppose $d_1>2$. We prove by induction on $1\leq j\leq d_{2}-d_{1}$ that:
\[J_{d_{1}+d_{3}-2+j} = Shad(J_{d_1+d_3-3+j})\cup \{x_{3}^{d_{3}+d_{1}-d_{2}-2+2j}x_{2}^{d_{2}-d_{1}+1-j}
\{x_{1},x_{2}\}^{d_1 - 1}\} = \]
\[ = \{\{x_{1},x_{2}\}^{d_{1}+d_{3}-2+j},\ldots, x_{3}^{d_{3}+d_{1}-d_{2}-2+2j}\{x_{1},x_{2}\}^{d_{2}-j},\]\[
x_{1}^{d_{1}}x_{3}^{d_{3}+d_{1}-d_{2}-1+2j}\{x_{1},x_{2}\}^{d_{2}-d_{1}-j-1},\ldots ,x_{1}^{d_{1}}x_{3}^{d_{3}-2+j}\}.\]
We already proved this for $j=1$. Suppose the assertion is true for some $j<d_2-d_1$. Since
$|J_{d_{1}+d_{3}-1+j}|-|Shad(J_{d_{1}+d_{3}-2+j})|=d_1$ we must add $d_1$ new generators to $Shad(J_{d_{1}+d_{3}-2+j})$
and these new generators are $x_{3}^{d_{3}+d_{1}-d_{2}-2+2j}x_{2}^{d_{2}-d_{1}+1-j}\{x_{1},x_{2}\}^{d_1-1}$ because $J$ is strongly stable and $x_3$ is a strong Lefschetz element for $S/J$. In particular, we get:
\[J_{d_{2}+d_{3}-2} = \{ \{x_{1},x_{2}\}^{d_{2}+d_{3}-2},x_{3}\{x_{1},x_{2}\}^{d_{2}+d_{3}-3},\ldots ,x_{3}^{d_{3}+d_{2}-d_{1}-2}\{x_{1},x_{2}\}^{d_{1}}\},\]
which is the same formula as in the case $d_1=2$.

We prove by induction on $1\leq j \leq d_1$ that \[ J_{d_2+d_{3}-2+j} = Shad(J_{d_2+d_{3}-3+j})\cup \{ x_{3}^{d_{3}+d_2-d_1+2j-2}\{x_{1},x_{2}\}^{d_1-j} \}.\]
If $j=1$ then $|J_{d_2+d_{3}-1}|- |Shad(J_{d_2+d_{3}-2})| = d_1$ so we must add $d_1$ generators, which are precisely the elements of the set $x_{3}^{d_{3}}\{x_{1},x_{2}\}^{d-1}$ since $J$ is strongly stable and $x_3$ is a strong Lefschetz element for $S/J$. The induction step is similar and finally we obtain that $J_{d_{3}+2d-2}=S_{d_{3}+2d-2}$ and thus we cannot add new minimal generators of $J$ in degrees $>d_{3}+2d-2$.
\end{proof}

\begin{cor}
In the conditions of the above proposition, the number of minimal generators of $J$ is $1+d_1 + d_1d_2$.
\end{cor}

\begin{ex}
Let $d_{1}=3$, $d_{2}=4$ and $d_{3}=9$. Then 
\[ J = (x_{1}^{3},\;\; x_{1}^{2}x_{2}^{2},\;\; x_{1}x_{2}^{4},\;\; x_{2}^{6},\;\; x_{2}^{5}x_3,\;\; x_{3}^{6}x_{2}^{3}\{x_{1},x_{2}\},\;\; \] 
\[ x_{3}^{8}x_{2}\{x_{1},x_{2}\}^{2},\;\; x_{3}^{10}\{x_{1},x_{2}\}^{2},\;\; x_{3}^{12}\{x_{1},x_{2}\},\;\; x_{3}^{14}).\]
\end{ex}

\section{\textbf{Case $d_1+d_2 > d_3+1$.}}

\medskip
\begin{itemize}
\item Subcase $d_1=d_2=d_3$.
\end{itemize} 
\medskip

\begin{prop}
\em Let $2\leq d:=d_1=d_2=d_3$ be positive integers. The Hilbert function of the standard graded complete intersection $A=K[x_{1},x_{2},x_{3}]/I$, where $I$ is the ideal generated by $f_{1}$, $f_{2}$, $f_{3}$, with $f_i$ homogeneous polynomials of degree $d_i$, for all $i$, with $1\leq i\leq 3$, has the form:  
\begin{enumerate}
	\item $H(A,k) = \binom{k+2}{k}$, for $k\leq d-1$.
	\item $H(A,k) = \binom{k+2}{2} - \frac{3j(j+1)}{2}$, for $k=j+d-1$, 
	      where $0\leq j \leq \left\lfloor \frac{d-1}{2} \right\rfloor$.	      
  \item $H(A,k) = H(A,3d-k-3)$, for $k\geq \left\lceil  \frac{3d-3}{2} \right\rceil$.
\end{enumerate}
\end{prop}

\begin{proof}
It follows from \cite[Lemma 2.9(b)]{HPV}.
\end{proof}

\begin{cor}
\em In the conditions of Proposition $3.1$, let $J=\Gin(I)$ be the generic initial ideal of $I$ with respect to the reverse lexicographic order. If we denote by $J_{k}$ the set of monomials of $J$ of degree $k$, then:
\begin{enumerate}
	\item $|J_{k}| = 0$, for $k\leq d-1$.
	\item $|J_{k}| = \frac{3j(j+1)}{2}$, for $k=d-1+j$, where $0 \leq j \leq [\frac{3d-1}{2}]$.
	\item If $d$ is even, then $|J_{k}| = \frac{3d^{2}+3d(4j+2)}{8} + \frac{3j(j+1))}{2}$, for
	      $k = j + \frac{3d-2}{2}$, where $0\leq j\leq \frac{d-2}{2}$.        
	      
	      \noindent
	      If $d$ is odd, then $|J_{k}| = \frac{3(d^{2}-1)+12jd}{8} + \frac{3j^{2}}{2}$, for
        $k = j + \frac{3d-3}{2}$, where $0\leq j\leq \frac{d-1}{2}$.        
  \item $|J_{k}| = \frac{3d(d-1)}{2} + 3jd$, for $k=j+2d-2$, where $0 \leq j \leq d-1$.
  \item $J_{k} = S_{k}$, for $k\geq 3d-2$.
\end{enumerate}
\end{cor}

\begin{prop}
\em Let $2\leq d:=d_1=d_2=d_3$ be positive integers. Let $f_{1},f_{2},f_{3}\in K[x_{1},x_{2},x_{3}]$ be a regular sequence of homogeneous polynomials of degrees $d_{1},d_{2},d_{3}$. If $I=(f_{1},f_{2},f_{3})$ , $J=\Gin(I)$, the generic initial ideal with respect to the reverse lexicographic order, and $S/I$ has (SLP), then:
\[ J = ( x_{1}^{d-2}\{x_{1},x_{2}\}^{2}, x_{1}^{d-2j-1}x_{2}^{3j+1}, x_{1}^{d-2j-2}x_{2}^{3j+2}\; for\; 1\leq j \leq  \frac{d-3}{2} , \; x_{2}^{\frac{3d-1}{2}}, x_{3}x_{2}^{\frac{3d-3}{3}}, \]\[
x_{3}^{2j+1}x_{1}^{2j}x_{2}^{\frac{3d-3}{2}-3j}, \ldots,  x_{3}^{2j+1}x_{2}^{\frac{3d-3}{2}-j} , 
1\leq j \leq \frac{d-3}{2} \;
 ,x_{3}^{d-2+2j}\{x_{1},x_{2}\}^{d-j}, 1\leq j\leq d) \]
if $d$ is odd, or
\[ J = ( x_{1}^{d-2}\{x_{1},x_{2}\}^{2}, x_{1}^{d-2j-1}x_{2}^{3j+1}, x_{1}^{d-2j-2}x_{2}^{3j+2}\; for\; 1\leq j \leq  \frac{d-4}{2}, \; x_{1}x_{2}^{\frac{3d-4}{2}}, x_{2}^{\frac{3d-2}{2}},\; \] \[
x_{3}^{2j}x_{1}^{2j-1}x_{2}^{\frac{3d}{2}-3j}, \ldots, x_{3}^{2j}x_{2}^{\frac{3d-2}{2}-j},  1\leq j\leq\frac{d-2}{2}
 ,x_{3}^{d-2+2j}\{x_{1},x_{2}\}^{d-j}, 1\leq j\leq d) \]
if $d$ is even.
\end{prop}

\begin{proof}
We have $|J_{d}|=3$, hence 
$J_{d}=\{ x_{1}^{d-2}\{x_{1},x_{2}\}^{2} \}$, since $J$ is strongly stable and $x_3$ is strong Lefschetz for $S/J$. Therefore:
\[ 
Shad(J_{d}) = \{x_{1}^{d-2}\{x_{1},x_{2}\}^{3} ,x_{1}^{d-2}x_{3}\{x_{1},x_{2}\}^{2}\}.
\]
Now we have four possibilities to analyze: $d=2$, $d=3$, $d=4$ and $d\geq 5$.

$\mathbf{d=2}$. Using the formulae from Corrolary $3.2$ we have $|J_{3}|-|Shad(J_2)|=2$ so there are two
generators to add to $Shad(J_{2})$ to obtain $J_{3}$. Since $J$ is strongly stable and $x_3$ is a strong Lefschetz element for $S/J$ these new generators are $x_{3}^{2}x_{1}$ and $x_{3}^{2}x_{2}$. Therefore
\[ J_3 = \{\{x_{1},x_{2}\}^{3} ,x_{3}\{x_{1},x_{2}\}^{2}, x_{3}^{2}\{x_{1},x_{2}\}\}.\]
Since $|J_4|-|Shad(J_3)|=1$ there is only one generator to add to $Shad(J_3)$ and this is precisely $x_{3}^{4}$.
It follows $J_4=S_4$ and thus we cannot add new minimal generators of $J$ in degree $>5$.

$\mathbf{d=3}$. We have $|J_{4}| - |Shad(J_3)|= 2$ so there are two generators to add to $Shad(J_{3})$ to obtain $J_{4}$. Since $J$ is strongly stable and $x_3$ is a strong Lefschetz element for $S/J$ these new generators are $x_{2}^{4}$ and $x_{3}x_{2}^{3}$. Therefore
\[ J_4 = \{\{x_{1},x_{2}\}^{4} ,x_{3}\{x_{1},x_{2}\}^{3}\}. \]
Since $|J_5|-|Shad(J_4)| = 3$ there are three new monomials to add to $Shad(J_4)$ in order to obtain $J_5$. Since
$J$ is strongly stable and $x_3$ is a strong Lefschetz element for $S/J$ these new generators are $x_{3}^{3}\{x_{1},x_{2}\}^{2}$. Analogously, we must add two new monomial to $Shad(J_5)$ in order to obtain $J_6$ and these are $x_{3}^{5}\{x_{1},x_{2}\}$. Finally, we will add $x_{3}^{7}$ and thus we cannot add new minimal generators of $J$ in degree $>7$.

$\mathbf{d=4}$. We have $|J_{4}| - |Shad(J_3)|= 2$ so there are two generators to add to $Shad(J_{4})$ to obtain $J_{5}$. Since $J$ is strongly stable and $x_3$ is a strong Lefschetz element for $S/J$ these new generators are $x_{1}x_{2}^{4}$ and $x_{2}^{5}$. Therefore
\[ J_5 = \{\{x_{1},x_{2}\}^{5} ,x_{1}^{2}x_{3}\{x_{1},x_{2}\}^{2}\}. \]
Since $|J_6|-|Shad(J_5)| = 2$ there are two new monomials to add to $Shad(J_5)$ in order to obtain $J_6$ and using the usual argument these new monomials are $x_{3}^{2}x_{1}x_{2}^{2},x_{3}^{2}x_{2}^{4}$. It follows
\[ J_6 = \{\{x_{1},x_{2}\}^{6} , x_{3}\{x_{1},x_{2}\}^{5},  x_{3}^{2}\{x_{1},x_{2}\}^{4}\}. \]
Finally, we will add consequently $x_{3}^{4}\{x_{1},x_{2}\}^{3}$, $x_{3}^{6}\{x_{1},x_{2}\}^{2}$, 
$x_{3}^{8}\{x_{1},x_{2}\}$ and $x_{3}^{10}$.

Suppose now $d\geq 5$. We have $|J_{d+1}|-|Shad(J_{d})|=2$ so there are two generators to add to $Shad(J_{d})$ to obtain $J_{d+1}$. Since $J$ is strongly stable and $x_3$ is a strong Lefschetz element for $S/J$ these new generators are
$x_{1}^{d-3}x_{2}^{4}, x_{1}^{d-4}x_{2}^{5}$. It follows
\[ J_{d+1} = \{ x_{1}^{d-4}\{x_{1},x_{2}\}^{5}, x_{1}^{d-2}x_{3}\{x_{1},x_{2}\}^{3}, x_{1}^{d-2}x_{3}^{2}\{x_{1},x_{2}\}^{2}\}. \]
We prove by induction on $j$, with $1\leq j \leq \left\lfloor \frac{d-3}{2} \right\rfloor$ that
\[ J_{d+j} = Shad(J_{d+j-1})\cup \{ x_{1}^{d-2j-1}x_{2}^{3j+1}, x_{1}^{d-2j-2}x_{2}^{3j+2} \} = \]
\[ = \{x_{1}^{d-2j-2}\{x_{1},x_{2}\}^{3j+2}, x_{1}^{d-2j+1}x_{3}\{x_{1},x_{2}\}^{3j-1}, \ldots, x_{1}^{d-2}x_{3}^{j}\{x_{1},x_{2}\}^{2}\}, \]
the assertion being checked for $j=1$. Assume now that the statement is true for 
some $j \leq \left\lfloor \frac{d-3}{2} \right\rfloor$. Then
\[Shad(J_{d+j}) = \{x_{1}^{d-2j-2}\{x_{1},x_{2}\}^{3j+3}, x_{1}^{d-2j-2}x_{3}\{x_{1},x_{2}\}^{3j+2},  \ldots, x_{1}^{d-2}x_{3}^{j+1}\{x_{1},x_{2}\}^{2}\}.\]
Since $|J_{d+j+1}|-|Shad(J_{d+j})|=2$ we must add two generators to $Shad(J_{d+j})$ to obtain $J_{d+j+1}$. Using the
fact that $J$ is strongly stable and $x_3$ is a strong Lefschetz element for $S/J$ it follows that these new generators are
$x_{1}^{d-2j-3}x_{2}^{3j+4}, x_{1}^{d-2j-4}x_{2}^{3j+5}$, so the induction step is fulfilled.

We must consider now two possibilities.

1. $d$ is odd. We obtain
	\[ J_{\frac{3d-3}{2}} = \{x_{1}\{x_{1},x_{2}\}^{\frac{3d-5}{2}}, x_{1}^{3}x_{3}\{x_{1},x_{2}\}^{\frac{3d-11}{2}}, 
	\ldots, x_{1}^{d-2}x_{3}^{\frac{d-3}{2}}\{x_{1},x_{2}\}^{2}\}.\]
	Since $|J_{\frac{3d-1}{2}}|-|Shad(J_{\frac{3d-3}{2}})|=2$ there are two generators to add to
	$Shad(J_{\frac{3d-3}{2}})$ to obtain $J_{\frac{3d-1}{2}}$, and they must be $x_{2}^{\frac{3d-1}{2}}$,
	$x_{3}x_{2}^{\frac{3d-3}{3}}$ using the usual argument. Therefore,
	\[ J_{\frac{3d-1}{2}} = \{ \{x_{1},x_{2}\}^{\frac{3d-1}{2}}, x_{3}\{x_{1},x_{2}\}^{\frac{3d-3}{2}}, 
	\ldots, x_{1}^{d-2}x_{3}^{\frac{d-1}{2}}\{x_{1},x_{2}\}^{2}\}.\]
  Since $|J_{\frac{3d+1}{2}}|-|Shad(J_{\frac{3d-1}{2}})|=1$ we must add a $3$ new generators to 
  $Shad(J_{\frac{3d-1}{2}})$ to obtain $J_{\frac{3d+1}{2}}$ and since $J$ is strongly stable and $x_3$ is a strong Lefschetz element for $S/J$, they are
  $x_{1}^{2}x_{3}^{3}x_{2}^{\frac{3d-9}{2}},x_{1}x_{3}^{3}x_{2}^{\frac{3d-7}{2}},x_{3}^{3}x_{2}^{\frac{3d-5}{2}}$.
  
We prove by induction on $j$, with $1\leq j \leq  \frac{d-3}{2}$ that
\[J_{\frac{3d-1}{2}+j} = \{Shad(J_{\frac{3d-3}{2}+j})\}\cup \{x_{3}^{2j+1}x_{1}^{2j}x_{2}^{\frac{3d-3}{2}-3j}, \ldots,  x_{3}^{2j+1}x_{2}^{\frac{3d-3}{2}-j}\} =\]\[ = \{ \{x_{1},x_{2}\}^{\frac{3d-1}{2}+j}, x_{3}\{x_{1},x_{2}\}^{\frac{3d-3}{2}+j}, \ldots, x_{3}^{2j+1}\{x_{1},x_{2}\}^{\frac{3d-3}{2}-j},\] 
\[ x_{3}^{2j+2}x_{1}^{2j+3}\{x_{1},x_{2}\}^{\frac{3d-11}{2}-3j},\ldots,x_{1}^{d-2}x_{3}^{\frac{d-1}{2}+j}\{x_{1},x_{2}\}^{2}\}. \]
This assertion is proved for $j=1$. Assume the assertion is true for some $j< \frac{d-3}{2} $.
Since $|J_{\frac{3d+1}{2}+j}| - |Shad(J_{\frac{3d-1}{2}+j})| = 2j+3$ we must add $2j+3$ new generators to 
$Shad(J_{\frac{3d-1}{2}+j})$ in order to obtain $J_{\frac{3d+1}{2}+j}$. The usual argument implies that those new generators are $x_{3}^{2j+3}x_{1}^{2j+2}x_{2}^{\frac{3d-3}{2}-3j-3}, \ldots, x_{3}^{2j+3}x_{2}^{\frac{3d-3}{2}-j-1}$, which conclude the induction.

2. $d$ is even. We obtain
\[ J_{\frac{3d-4}{2}} = \{x_{1}^{2}\{x_{1},x_{2}\}^{\frac{3d-8}{2}}, x_{1}^{4}x_{3}\{x_{1},x_{2}\}^{\frac{3d-14}{2}}, 
\ldots, x_{1}^{d-2}x_{3}^{\frac{d-4}{2}}\{x_{1},x_{2}\}^{2}\}.\]
We have $|J_{\frac{3d-2}{2}}|-|Shad(J_{\frac{3d-4}{2}})|=2$ so we must add two new generators to 
$Shad(J_{\frac{3d-4}{2}})$ to obtain $J_{\frac{3d-2}{2}}$. Since $J$ is strongly stable and $x_3$ is a strong Lefschetz element for $S/J$, they are $x_{1}x_{2}^{\frac{3d-4}{2}}$ and $x_{2}^{\frac{3d-2}{2}}$, therefore
\[ J_{\frac{3d-2}{2}} = \{\{x_{1},x_{2}\}^{\frac{3d-2}{2}}, x_{1}^{2}x_{3}\{x_{1},x_{2}\}^{\frac{3d-8}{2}}, \ldots, x_{1}^{d-2}x_{3}^{\frac{d-2}{2}}\{x_{1},x_{2}\}^{2}\}.\]
Since $|J_{\frac{3d}{2}}|-|Shad(J_{\frac{3d-2}{2}})|=2$ we must add two new generators to $Shad(J_{\frac{3d-2}{2}})$
in order to obtain $J_{\frac{3d}{2}}$ and since $J$ is strongly stable and $x_3$ is a strong Lefschetz element for $S/J$, they are $x_{1}x_{3}^{2}x_{2}^{\frac{3d-6}{2}}$ and
$x_{3}^{2}x_{2}^{\frac{3d-4}{2}}$.

We prove by induction on $j$, with $1\leq j \leq  \frac{d-2}{2} $ that
\[ J_{\frac{3d-2}{2} + j} = Shad(J_{\frac{3d-4}{2}+j}) \cup \{ x_{3}^{2j}x_{1}^{2j-1}x_{2}^{\frac{3d}{2}-3j}, \ldots, x_{3}^{2j}x_{2}^{\frac{3d-2}{2}-j}\}.\]
The assertion has been proved for $j=1$ and suppose it is true for some \linebreak $j< \frac{d-2}{2} $.
Since $|J_{\frac{3d}{2} + j}| - |Shad(J_{\frac{3d-2}{2} + j})| = 2j+2$ we must add $2j+2$ generators to
$Shad(J_{\frac{3d-2}{2} + j})$ in order to obtain $J_{\frac{3d}{2} + j}$ and since $J$ is strongly stable and $x_3$ is a strong Lefschetz element for $S/J$ they must be 
$x_{3}^{2j+2}x_{1}^{2j+1}x_{2}^{\frac{3d}{2}-3j-3}, \ldots, x_{3}^{2j+2}x_{2}^{\frac{3d-4}{2}-j}$, which conclude
the induction.

Either if $d$ is even, either if $d$ is odd, we obtain
\[J_{2d-2} = \{\{x_{1},x_{2}\}^{2d-2},x_{3}\{x_{1},x_{2}\}^{2d-2},\ldots,x_{3}^{d-2}\{x_{1},x_{2}\}^{d}\} = \{\{x_{1},x_{2}\}^{d} \{x_{1},x_{2},x_{3}\}^{d-2}\}.\]
Since $|J_{2d-1}|-|Shad(J_{2d-2})|=d$ we must add $d$ new generators to $Shad(J_{2d-2})$ to obtain $J_{2d-1}$. But
$J$ is strongly stable and $x_3$ is a strong Lefschetz element for $S/J$, so we must add $x_{3}^{d}\{x_{1},x_{2}\}^{d-1}$, therefore $J_{2d-1}=\{\{x_{1},x_{2}\}^{d-1} \{x_{1},x_{2},x_{3}\}^{d}\}$. Using induction on $j\leq d$, we prove that
 \[ J_{2d-2+j} = Shad(J_{2d-3+j}) \cup \{x_{3}^{d-2+2j}\{x_{1},x_{2}\}^{d-j}\} = \{x_{1},x_{2}\}^{d-j}\{x_{1},x_{2},x_{3}\}^{d-2+2j}.\]
For $j=1$ we already proved. Suppose the assertion is true for $j<d$. We have $|J_{2d-1+j}|-|Shad(J_{2d-2+j})|=d-j$ so
we must add $d-j$ new monomials to $Shad(J_{2d-2+j})$ to obtain $J_{2d-1+j}$ and from the usual argument, these new monomials are $x_{3}^{d+2j}\{x_{1},x_{2}\}^{d-j-1}$. Finally, since $J_{3d-2}=S_{3d-2}$ we cannot add new minimal generators of $J$ in degree $>3d-2$.
\end{proof}

\begin{cor}
In the conditions above, the number of minimal generators of $J$ is 
$1 + \frac{d(d+1)}{2} + \frac{(d+1)^{2}}{4}$ when $d$ is odd, 
or $1 + \frac{d(d+1)}{2} + \frac{d(d+2)}{4}$ when $d$ is even.
\end{cor}

\begin{ex}
\begin{enumerate}
	\item Let $d_{1}=d_{2}=d_{3}=5$. Then
	\[ J = (x_{1}^{3}\{x_1,x_2\}^{2},\; x_1^{2}x_2^{4},\;x_1x_2^{5},\; x_2^{7},\; x_3 x_2^{6},\;
	        x_{2}^{3}x_{3}^{3}\{x_1,x_2\}^{2},\; \]
	\[ x_{3}^{5}\{x_1,x_2\}^{4},\;x_{3}^{7}\{x_1,x_2\}^{3},\; x_{3}^{9}\{x_1,x_2\}^{2}, \; x_{3}^{11}\{x_1,x_2\}, \;
	x_{3}^{13}).\]
	\item Let $d_{1}=d_{2}=d_{3}=6$. Then
	\[ J = (x_1^{4}\{x_{1},x_{2}\}^{2},\; x_{1}^{3}x_{2}^{4},\; x_{1}^{2}x_{2}^{5},\; x_{1}x_{2}^{7},\; x_{2}^{8},\;
	        x_2^{6}x_3^{2}\{x_1,x_2\},\;x_2^{3}x_3^{4}\{x_1,x_2\}^{3},\;\]
  \[ x_{3}^{6}\{x_{1},x_2\}^{5},\; x_{3}^{8}\{x_{1},x_2\}^{4},\; x_{3}^{10}\{x_{1},x_2\}^{3},\;
	         x_{12}^{6}\{x_{1},x_2\}^{2},\; x_{14}^{6}\{x_{1},x_2\},\; x_{16}^{6}) \]
\end{enumerate}
\end{ex}

\medskip
\begin{itemize}
\item Subcase $d_1=d_2<d_3$.
\end{itemize} 
\medskip

\begin{prop}
\em Let $2\leq d:=d_1=d_2<d_3$ be positive integers such that $d_1+d_2 > d_3 + 1$. The Hilbert function of the standard graded complete intersection $A=K[x_{1},x_{2},x_{3}]/I$, where $I$ is the ideal generated by $f_{1}$, $f_{2}$, $f_{3}$, with $f_i$ homogeneous polynomials of degree $d_i$, for all $i$, with $1\leq i\leq 3$, has the form:  
\begin{enumerate}
	\item $H(A,k) = \binom{k+2}{k}$, for $k\leq d-1$.
	\item $H(A,k) = \binom{d+1}{2} + \sum_{i=1}^{j}(d-i)$, for $k=j+d-1$, where $0\leq j \leq d_{3}-d$.
	\item $H(A,k) = \binom{d+1}{2} + \sum_{i=1}^{d_{3}-d}(d-i) + \sum_{i=1}^{j}(2d-d_{3}-2i)$, 
	      for $k=j + d_3 - 1$, where $0\leq j \leq \left\lfloor \frac{2d-d_{3}-1}{2} \right\rfloor$.
	\item $H(A,k) = H(A,d_{3}+2d-3-k)$,for $k\geq \left\lceil \frac{d_{3}+2d-3}{2} \right\rceil$.
\end{enumerate}
\end{prop}

\begin{proof}
It follows from \cite[Lemma 2.9(b)]{HPV}.
\end{proof}

\begin{cor}
\em In the conditions of Proposition $3.6$, let $J=\Gin(I)$ be the generic initial ideal of $I$ with respect to the reverse lexicographic order. If we denote by $J_{k}$ the set of monomials of $J$ of degree $k$, then:
\begin{enumerate}
	\item $|J_{k}| = 0$, for $k\leq d-1$.
	\item $|J_{k}| = j(j+1)$, for $k = j + d - 1$ where $0 \leq j \leq d_{3}-d$.
	\item $|J_{k}| = d_{3}^{2}+d_{3}-d^{2}-d-2dd_{3} + j(2d_{3}-d) + \frac{3j(j+1)}{2}$,
	      for $k=j+d_{3}-1$, where $0\leq j \leq \left\lfloor \frac{2d-d_{3}-1}{2} \right\rfloor$.
	
	\item If $d_{3}$ is even then $|J_{k}| = \frac{4d^{2}+3d_{3}^{2}-4dd_{3}+4d}{8} + \frac{j(2d+d_{3})}{2} +
	      \frac{3j(j+1)}{2}$, for $k = j + \frac{2d+d_{3}-2}{2}$, where $0 \leq j \leq \frac{2d-d_{3}-2}{2}$.
	      	
	      If $d_{3}$ is odd then $|J_{k}| = \frac{3d_{3}^{2}+4d^{2}-4dd_{3}-3}{2} + \frac{j(2d+d_{3}-3)}{2} +
	      \frac{3j(j+1)}{2}$, for $k = j + \frac{2d+d_{3}-3}{2}$, where $0 \leq j \leq \frac{2d-d_{3}-1}{2}$.	      		
  \item $|J_{k}| =  3d^{2} - 2d + \frac{d_{3}(d_{3}+1)}{2} - 2dd_{3} + (4d-d_{3})j + j(j-1)$, for $k=j+2d-2$,
        where $0 \leq j \leq d_{3}-d$.        
  \item $|J_{k}| = \frac{(d+d_{3})(d+d_{3}-1)}{2} - \frac{d(d+1)}{2} + j(2d+d_{3})$, for $k=j+d_{3}+d-2$,
        where $0 \leq j \leq d-1$.	
  \item $J_{k} = S_{k}$, for $k\geq 3d-2$.
\end{enumerate}
\end{cor}

\begin{prop}
\em Let $2\leq d:=d_1=d_2<d_3$ be positive integers such that $2d>d_3+1$. Let $f_{1},f_{2},f_{3}\in K[x_{1},x_{2},x_{3}]$ be a regular sequence of homogeneous polynomials of degrees $d_{1},d_{2},d_{3}$. If $I=(f_{1},f_{2},f_{3})$, $J=\Gin(I)$, the generic initial ideal with respect to the reverse lexicographic order and
$S/I$ has (SLP), then if $d_3$ is even, we have:
	\[ J = (x_{1}^{d},x_{1}^{d-1}x_{2}, x_{1}^{d-j-1}x_{2}^{2j+1} \;for\; 1\leq j\leq d_{3}-d-1, \]
	\[ x_{1}^{2d-d_{3}-2j+1}x_{2}^{2d_{3}-2d+3j-2}, x_{1}^{2d-d_{3}-2j}x_{2}^{2d_{3}-2d+3j-1}\;for\;1\leq j\leq
	    \frac{2d-d_{3}}{2} , \]
	\[ x_{3}^{2j}x_{1}^{2j-1}x_{2}^{\frac{2d+d_{3}-2}{2}-3j}, \ldots, x_{3}^{2j}x_{2}^{\frac{2d+d_{3}-4}{2}-j} 
	\;for\;1\leq j\leq \frac{2d-d_{3}-2}{2}, \]
	\[ x_{3}^{2d-d_{3}-2+2j}x_{2}^{2d_{3}-2d+2-2j}\{x_{1},x_{2}\}^{2d-d_{3}+j-2},1\leq j\leq d_{3}-d, 
	   x_{3}^{d_{3}-2+2j}\{x_{1},x_{2}\}^{d-j},1\leq j\leq d).\]
	Otherwise, if $d_3$ is odd, we have
	\[ J = (x_{1}^{d},x_{1}^{d-1}x_{2}, x_{1}^{d-j-1}x_{2}^{2j+1} \;for\; 1\leq j\leq d_{3}-d-1, \]
	\[ x_{1}^{2d-d_{3}-2j+1}x_{2}^{2d_{3}-2d+3j-2}, x_{1}^{2d-d_{3}-2j}x_{2}^{2d_{3}-2d+3j-1}\;for\;1\leq j\leq
	    \frac{2d-d_{3}-1}{2} , \]
	\[ x_{2}^{\frac{2d+d_{3}-1}{2}}, x_{3}x_{2}^{\frac{2d+d_{3}-3}{2}}, x_{3}^{2j+1}x_{1}^{2j}x_{2}^{\frac{2d+d_{3}-3}{2}-3j}, \ldots, x_{3}^{2j+1}x_{2}^{\frac{2d+d_{3}-3}{2}-j}
	,1\leq j\leq \frac{2d-d_{3}-3}{2}, \]
	\[ x_{3}^{2d-d_{3}-2+2j}x_{2}^{2d_{3}-2d+2-2j}\{x_{1},x_{2}\}^{2d-d_{3}+j-2},1\leq j\leq d_{3}-d, x_{3}^{d_{3}-2+2j}\{x_{1},x_{2}\}^{d-j},1\leq j\leq d). \]	
\end{prop}

\begin{proof}
We note first that $d\geq 3$. Indeed, if $d=2$ then the condition $2d>d_{3}+1$ implies $d_{3}=2$ which is a contradiction. We have $|J_{d}|=2$, hence $J_{d}=x_{1}^{d-1}\{x_{1},x_{2}\}$, since $J$ is strongly stable. Therefore:
\[ Shad(J_{d}) = \{x_{1}^{d-1}\{x_{1},x_{2}\}^{2} ,x_{1}^{d-1}x_{3}\{x_{1},x_{2}\}\}. \]
Assume $d_3>d+1$. Since $|J_{d+1}| - |Shad(J_{d})| = 1$ we must add a new generator to $Shad(J_{d})$ in order to obtain $J_{d+1}$. On the other hand, $J$ is strongly stable and $x_3$ is a strong Lefschetz element for $S/J$ so this new generator is $x_{1}^{d-2}x_{2}^{3}$, therefore
\[ J_{d+1} = \{x_{1}^{d-2}\{x_{1},x_{2}\}^{3}, x_{1}^{d-1}x_{3}\{x_{1},x_{2}\}\}.\]
We prove by induction on $1\leq j\leq d_{3}-d-1$ that
\[J_{d+j} = Shad(J_{d-1+j})\cup \{x_{1}^{d-j-1}x_{2}^{2j+1}\} = \{x_{1}^{d-j-1}\{x_{1},x_{2}\}^{2j+1}, \ldots, x_{1}^{d-1}x_{3}^{j}\{x_{1},x_{2}\}\}.\]
The case $j=1$ was done. Suppose the assertion is true for some $j<d_{3}-d-1$. Then, since $|J_{d+j+1}|-|Shad(J_{d+j})|=1$
it follows that we must add one generator to $Shad(J_{d+j})$ in order to obtain $J_{d+j+1}$. Since $J$ is strongly stable and $x_3$ is a strong Lefschetz element for $S/J$,
this new generator must be $x_{1}^{d-j-2}x_{2}^{2j+3}$. In particular,
\[ J_{d_{3}-1} = \{x_{1}^{2d-d_{3}}\{x_{1},x_{2}\}^{2d_{3}-2d-1}, x_{1}^{2d-d_{3}+1}x_{3}\{x_{1},x_{2}\}^{2d_{3}-2d-3}, \ldots, x_{3}^{d_{3}-d-1}x_{1}^{d-1}\{x_{1},x_{2}\}\}, \]
which is the same formula as in the case $d_3=d+1$.

We need to consider several possibilities. First, suppose $d_3 = 2d-2$. We have $|J_{d_{3}}| - |Shad(J_{d_{3}-1})| = 2$ so we must add two generators to $Shad(J_{d_{3}-1})$ in order to obtain $J_{d_{3}} = J_{2d-2}$. But since $J$ is strongly
stable and $x_3$ is a strong Lefschetz element for $S/J$, these new generators are $x_{1}x_{2}^{2d-3}$ and $x_{2}^{2d-2}$.

Suppose now $d_3<2d-2$. Since $|J_{d_{3}}| - |Shad(J_{d_{3}-1})| = 2$ we must add two generators to $Shad(J_{d_{3}-1})$
in order to obtain $J_{d_{3}}$. $J$ strongly stable and $x_3$ is a strong Lefschetz element for $S/J$ force us to choose $x_{1}^{2d-d_{3}-1}x_{2}^{2d_{3}-2d+1}, x_{1}^{2d-d_{3}-2}x_{2}^{2d_{3}-2d+2}$, so
\[ J_{d_{3}} = \{x_{1}^{2d-d_{3}-2}\{x_{1},x_{2}\}^{2d_{3}-2d}, x_{1}^{2d-d_{3}}x_{3}\{x_{1},x_{2}\}^{2d_{3}-2d-1}, x_{1}^{2d-d_{3}+1}x_{3}^{2}\{x_{1},x_{2}\}^{2d_{3}-2d-3}\}. \]
We show by induction on $1\leq j\leq \left\lfloor \frac{2d-d_{3}+1}{2} \right\rfloor$ that
\[ J_{d_{3}-1+j} = Shad(J_{d_{3}-2+j})\cup \{ x_{1}^{2d-d_{3}-2j+1}x_{2}^{2d_{3}-2d+3j-2}, x_{1}^{2d-d_{3}-2j}x_{2}^{2d_{3}-2d+3j-1}\} = \]
\[ = \{ x_{1}^{2d-d_{3}-2j}\{x_{1},x_{2}\}^{2d_{3}-2d+3j-1}, x_{3}x_{1}^{2d-d_{3}-2j+2}\{x_{1},x_{2}\}^{2d_{3}-2d+3j-4}, \ldots,\]\[ x_{3}^{j}x_{1}^{2d-d_{3}}\{x_{1},x_{2}\}^{2d_{3}-2d-1},
x_{3}^{j+1}x_{1}^{2d-d_{3}+1}\{x_{1},x_{2}\}^{2d_{3}-2d-3},\ldots, x_{3}^{d_{3}-d-1+j}x_{1}^{d-1}\{x_{1},x_{2}\}\}.\]
We already done the case $j=1$. Suppose the assertion is true for some \linebreak $j<\left\lfloor \frac{2d-d_{3}+1}{2} \right\rfloor$. We have $|J_{d_{3}+j}|-|Shad(J_{d_{3}+j-1})|=2$ so we add two generators to $Shad(J_{d_{3}+j-1})$ and they must be $x_{1}^{2d-d_{3}-2j-1}x_{2}^{2d_{3}-2d+3j+1}, x_{1}^{2d-d_{3}-2j-2}x_{2}^{2d_{3}-2d+3j+2}$
from the usual argument.
In the following, we distinguish between two possibilities: $d$ is even or $d$ is odd. If $d$ is even, we get
\[ J_{\frac{2d+d_{3}-4}{2}} = \{ x_{1}^{2}\{x_{1},x_{2}\}^{\frac{2d+d_{3}-8}{2}}, x_{3}x_{1}^{4}\{x_{1},x_{2}\}^{\frac{2d+d_{3}-14}{2}}, \ldots, \]
\[ x_{3}^{\frac{2d-d_{3}-2}{2}}x_{1}^{2d-d_{3}}\{x_{1},x_{2}\}^{2d_{3}-2d-1},\ldots
,x_{3}^{\frac{d_{3}-3}{2}}x_{1}^{d-1}\{x_{2},x_{2}\} \}.\]
We have $|J_{\frac{2d+d_{3}-2}{2}}|-|Shad(J_{\frac{2d+d_{3}-4}{2}})|=2$ so we add two generators to $Shad(J_{\frac{2d+d_{3}-4}{2}})$ and they must be $x_{1}x_{2}^{\frac{2d+d_{3}-4}{2}}, x_{2}^{\frac{2d+d_{3}-2}{2}}$, so:
\[ J_{\frac{2d+d_{3}-2}{2}} =  \{\{x_{1},x_{2}\}^{\frac{2d+d_{3}-2}{2}}, x_{3}x_{1}^{2}\{x_{1},x_{2}\}^{\frac{2d+d_{3}-8}{2}}, \ldots, \]
\[,x_{3}^{\frac{2d-d_{3}}{2}}x_{1}^{2d-d_{3}}\{x_{1},x_{2}\}^{2d_{3}-2d-1},\ldots
,x_{3}^{\frac{d_{3}-2}{2}}x_{1}^{d-1}\{x_{2},x_{2}\}\}.\]
One can easily show by induction on  $1\leq j\leq \frac{2d-d_{3}-2}{2}$, if case, that
\[ J_{\frac{2d+d_{3}-2}{2}+j} = Shad(J_{\frac{2d+d_{3}-4}{2}+j})\cup \{ x_{3}^{2j}x_{1}^{2j-1}x_{2}^{\frac{2d+d_{3}-2}{2}-3j}, \ldots, x_{3}^{2j}x_{2}^{\frac{2d+d_{3}-4}{2}-j}\} = \]
\[ = \{ \{x_{1},x_{2}\}^{\frac{2d+d_{3}-2}{2}+j}, \ldots, x_{3}^{2j}\{x_{1},x_{2}\}^{\frac{2d+d_{3}-2}{2}-j},x_{3}^{2j+1}x_{1}^{2j+2}\{x_{1},x_{2}\}^{\frac{2d+d_{3}-8}{2}-3j},
 \ldots, \]\[
x_{3}^{\frac{2d-d_{3}}{2} +j}x_{1}^{2d-d_{3}}\{x_{1},x_{2}\}^{2d_{3}-2d-1}, \ldots, 
x_{3}^{\frac{d_{3}-2}{2} + j}x_{1}^{d-1}\{x_{1},x_{2}\} \}. \]
Indeed, if $j=1$, we have $|J_{\frac{2d+d_{3}}{2}+j}| - |Shad(J_{\frac{2d+d_{3}-2}{2}+j})| = 2$ so we must add two
monomials to $Shad(J_{\frac{2d+d_{3}-2}{2}+j})$ in order to obtain $J_{\frac{2d+d_{3}}{2}+j}$. Since $J$ is strongly stable and $x_3$ is a strong Lefschetz element for $S/J$, these new monomials are $x_{3}^{2}x_{1}x_{2}^{\frac{2d+d_{3}-2}{2}-3}$ and $x_{3}^{2}x_{2}^{\frac{2d+d_{3}-2}{2}-2}$. The induction step is similar.

If $d$ is odd, we have $ J_{\frac{2d+d_{3}-3}{2}} = \{x_{1}\{x_{1},x_{2}\}^{\frac{2d+d_{3}-5}{2}}, x_{3}x_{1}^{3}\{x_{1},x_{2}\}^{\frac{2d+d_{3}-11}{2}}, \ldots, $

$x_{3}^{\frac{2d-d_{3}-1}{2}}x_{1}^{2d-d_{3}}\{x_{1},x_{2}\}^{2d_{3}-2d-1},\ldots
,x_{3}^{\frac{d_{3}-3}{2}}x_{1}^{d-1}\{x_{2},x_{2}\} \}$.

Since $|J_{\frac{2d+d_{3}-1}{2}}|-|Shad(J_{\frac{2d+d_{3}-3}{2}})| = 2$ we add two generators to $Shad(J_{\frac{2d+d_{3}-3}{2}})$ in order to obtain $J_{\frac{2d+d_{3}-1}{2}}$. Since $J$ is strongly stable and $x_3$ is a strong Lefschetz element for $S/J$, these new monomials are $x_{2}^{\frac{2d+d_{3}-1}{2}}, x_{3}x_{2}^{\frac{2d+d_{3}-3}{2}}$, therefore:
\[ J_{\frac{2d+d_{3}-1}{2}} =  \{ \{x_{1},x_{2}\}^{\frac{2d+d_{3}-1}{2}}, x_{3}\{x_{1},x_{2}\}^{\frac{2d+d_{3}-3}{2}}, x_{3}^{2}x_{1}^{3}\{x_{1},x_{2}\}^{\frac{2d+d_{3}-11}{2}}, \ldots, \]
\[x_{3}^{\frac{2d-d_{3}+1}{2}}x_{1}^{2d-d_{3}}\{x_{1},x_{2}\}^{2d_{3}-2d-1},\ldots
,x_{3}^{\frac{d_{3}-1}{2}}x_{1}^{d-1}\{x_{2},x_{2}\} \}.\]
Assume $d_3 < 2d-3$ (otherwise $\frac{2d+d_3-1}{2}=2d-2$). $|J_{\frac{2d+d_{3}+1}{2}}| - |Shad(J_{\frac{2d+d_{3}-3}{2}})| = 3$ so we must add $3$ generators to $Shad(J_{\frac{2d+d_{3}-3}{2}})$ to 
obtain $J_{\frac{2d+d_{3}+1}{2}}$. The usual argument implies that they are
$ x_{1}^{2}x_{3}^{3}x_{2}^{\frac{2d+d_{3}-9}{2}}, x_{1}x_{3}^{3}x_{2}^{\frac{2d+d_{3}-7}{2}}, x_{3}^{3}x_{2}^{\frac{2d+d_{3}-5}{2}}$ and thus
\[ J_{\frac{2d+d_{3}+1}{2}} = \{\{x_{1},x_{2}\}^{\frac{2d+d_{3}+1}{2}}, x_{3}\{x_{1},x_{2}\}^{\frac{2d+d_{3}-1}{2}}, x_{3}^{2}\{x_{1},x_{2}\}^{\frac{2d+d_{3}-3}{2}}, x_{3}^{3}\{x_{1},x_{2}\}^{\frac{2d+d_{3}-5}{2}}, \ldots, \]
\[x_{3}^{\frac{2d-d_{3}+3}{2}}x_{1}^{2d-d_{3}}\{x_{1},x_{2}\}^{2d_{3}-2d-1},\ldots
,x_{3}^{\frac{d_{3}+1}{2}}x_{1}^{d-1}\{x_{2},x_{2}\}\}.\]

One can easily prove by induction on $1\leq j\leq \frac{2d-d_{3}-3}{2}$ that
\[ J_{\frac{2d+d_{3}-1}{2}+j} = Shad(J_{\frac{2d+d_{3}-3}{2}+j})\cup \{ x_{3}^{2j+1}x_{1}^{2j}x_{2}^{\frac{2d+d_{3}-3}{2}-3j}, \ldots, x_{3}^{2j+1}x_{2}^{\frac{2d+d_{3}-3}{2}-j} \} = \]
\[ = \{ \{x_{1},x_{2}\}^{\frac{2d+d_{3}-1}{2}+j},\ldots, x_{3}^{2j+1}\{x_{1},x_{2}\}^{\frac{2d+d_{3}-3}{2}-j},
x_{3}^{2j+2}x_{1}^{2j+3}\{x_{1},x_{2}\}^{\frac{2d+d_{3}-11}{2}-3j}, \]\[
x_{3}^{\frac{2d-d_{3}+1}{2} +j}x_{1}^{2d-d_{3}}\{x_{1},x_{2}\}^{2d_{3}-2d-1}, \ldots, 
x_{3}^{\frac{d_{3}-1}{2} + j}x_{1}^{d-1}\{x_{1},x_{2}\} \}. \]

In all cases above, we get: $J_{2d-2} = \{\{x_{1},x_{2}\}^{2d-2},x_{3}\{x_{1},x_{2}\}^{2d-3},\ldots,$ \linebreak
$x_{3}^{2d-d_{3}-2}\{x_{1},x_{2}\}^{d_{3}}
,x_{3}^{2d-d_{3}-1}x_{1}^{2d-d_{3}}\{x_{1},x_{2}\}^{2d_{3}-2d-1},\ldots,x_{3}^{d-2}x_{1}^{d-1}\{x_{1},x_{2}\}\}$.  

We have $|J_{2d-1}| - |Shad(J_{2d-2})| = 2d-d_{3}$ so we must add $2d-d_{3}$ new generators to $Shad(J_{2d-2})$ to
obtain $J_{2d-1}$. We get
$ J_{2d-1} = \{ \{x_{1},x_{2}\}^{2d-1}, x_{3}\{x_{1},x_{2}\}^{2d-2}, \ldots,$ \linebreak
$ x_{3}^{2d-d_{3}}\{x_{1},x_{2}\}^{d_{3}-1}
,x_{3}^{2d-d_{3}+1}x_{1}^{2d-d_{3}+1}\{x_{1},x_{2}\}^{2d_{3}-2d-3},\ldots,x_{3}^{d-1}x_{1}^{d-1}\{x_{1},x_{2}\}\}$.

One can easily show by induction of $1\leq j\leq d_{3}-d$ that
\[ J_{2d-2+j} = Shad(J_{2d-3+j})\cup \{ x_{3}^{2d-d_{3}-2+2j}x_{2}^{2d_{3}-2d+2-2j}\{x_{1},x_{2}\}^{2d-d_{3}+j-2}\} = \]
\[ = \{ \{x_{1},x_{2}\}^{2d-2+j},x_{3}\{x_{1},x_{2}\}^{2d-3+j},\ldots,x_{3}^{2d-d_{3}-2+2j}\{x_{1},x_{2}\}^{d_{3}-j},\]
\[x_{3}^{2d-d_{3}-1+2j}x_{1}^{2d-d_{3}+j}\{x_{1},x_{2}\}^{2d_{3}-2d-1-2j},\ldots,x_{3}^{d-2+j}x_{1}^{d-1}\{x_{1},x_{2}\}\}.\]
In particular, we get: $J_{d+d_{3}-2} = \{\{x_{1},x_{2}\}^{d+d_{3}-2},x_{3}\{x_{1},x_{2}\}^{d+d_{3}-3},\ldots ,x_{3}^{d_{3}-2}\{x_{1},x_{2}\}^{d}\}$. 
We have $|J_{d+d_3-1}|-|Shad(J_{d+d_3-2})|=d$ so we must add $d$ new generators to $Shad(J_{d+d_3-2})$ in order
to obtain $J_{d+d_3-1}$. Since $J$ is strongly stable and $x_3$ is a strong Lefschetz element for $S/J$, these new generators are $x_{3}^{d_{3}}\{x_{1},x_{2}\}^{d-1}$ so
\[J_{d+d_{3}-1} = \{\{x_{1},x_{2}\}^{d+d_{3}-1},x_{3}\{x_{1},x_{2}\}^{d+d_{3}-2},\ldots ,x_{3}^{d_{3}}\{x_{1},x_{2}\}^{d-1}\}. \]
Now, one can easily prove by induction on $1\leq j\leq d-1$ that $J_{d+d_{3}-2+j}$ is the set
\[ Shad(J_{d+d_{3}-3+j}) \cup \{x_{3}^{d_{3}-2+2j}\{x_{1},x_{2}\}^{d-j}\} =
 \{ \{x_{1},x_{2}\}^{d+d_{3}-2+j},\ldots ,x_{3}^{d_{3}-2+2j}\{x_{1},x_{2}\}^{d-j}\}.\]
Finally we obtain that $J_{d_{3}+2d-2}=S_{d_{3}+2d-2}$ and thus we cannot add new minimal generators of $J$ in degree $>d_{3}+2d-2$.
\end{proof}

\begin{cor}
In the conditions of the above proposition, the number of minimal generators of $J$ is $d(d+1) - \left(\frac{2d-d_3}{2}\right)^{2} + 1$ if $d_3$ is even or $d(d+1) - \frac{(2d-d_3)^{2}-1}{4} + 1$ if $d_3$ is odd.
\end{cor}

\begin{ex}
\begin{enumerate}
	\item Let $d_{1}=d_{2}=4$ and $d_{3}=6$. We have
	\[J = (x_{1}^{4},\; x_{1}^{3}x_2,\; x_1^{2}x_2^{3},\; x_1 x_2^{5},\;x_2^{6},\; x_3^{2}x_2^{4}\{x_1,x_2\},\; 
	       x_{3}^{4}x_{2}^{2}\{x_1,x_2\}^{2},\]\[\; x_3^{6}\{x_1,x_2\}^{3},\; x_3^{8}\{x_1,x_2\}^{2},\;
	       x_3^{10}\{x_1,x_2\},\; x_3^{12}). \]
	\item Let $d_{1}=d_{2}=4$ and $d_{3}=5$. We have:
	\[J = (x_{1}^{4},\; x_{1}^{3}x_2,\; x_1^{2}x_2^{3},\; x_1 x_2^{4},\; x_2^{6},\;x_2^{5}x_{3},\; 
	       x_3^{3}x_2^{2}\{x_1,x_2\}^{2},\;\]\[\; x_3^{5}\{x_1,x_2\}^{3},\; x_3^{7}\{x_1,x_2\}^{2},\;
	       x_3^{9}\{x_1,x_2\},\; x_3^{11} ). \]
\end{enumerate}
\end{ex}

\medskip
\begin{itemize}
\item Subcase $d_1<d_2=d_3$.
\end{itemize} 
\medskip

\begin{prop}
\em Let $2\leq d_1<d_2=d_3=:d$ be positive integers. The Hilbert function of the standard graded complete intersection $A=K[x_{1},x_{2},x_{3}]/I$, where $I$ is the ideal generated by $f_{1}$, $f_{2}$, $f_{3}$, with $f_i$ homogeneous polynomials of degree $d_i$, for all $i$, with $1\leq i\leq 3$, has the form:  
\begin{enumerate}
	\item $H(A,k) = \binom{k+2}{k}$, for $k\leq d_{1}-2$.
	\item $H(A,k) = \binom{d_{1}+1}{2} + jd_{1}$, for $k = j +d_{1} - $, where $0 \leq j \leq d-d_1$.
	\item $H(A,k) = \binom{d_{1}+1}{2} + d_{1}(d-d_{1}) + \sum_{i=1}^{j}(d_{1}-2i)$, 
        for $k = j + d - 1$, where $0 \leq j \leq \left\lfloor \frac{d_{1}-1}{2} \right\rfloor$.
	\item $H(A,k) = H(A,d_{1}+2d-3-k)$, for $k\geq \left\lceil \frac{d_{1}+2d-3}{2} \right\rceil$.
\end{enumerate}
\end{prop}

\begin{proof}
It follows from \cite[Lemma 2.9(b)]{HPV}.
\end{proof}

\begin{cor}
\em In the conditions of Proposition $3.11$, let $J=\Gin(I)$ be the generic initial ideal of $I$ with respect to the reverse lexicographic order. If we denote by $J_{k}$ the set of monomials of $J$ of degree $k$, then:
\begin{enumerate}
	\item $|J_{k}| = 0$, for $k\leq d_{1}-1$.
	\item $|J_{k}| = j(j+1)/2$, for $k=j+d_{1}-1$, where $0 \leq j \leq d-d_{1}$.	
	\item $|J_{k}| = \frac{(d-d_{1})(d-d_{1}-1)}{2} + j(d-d_{1}) + \frac{3j(j+1)}{2}$,
	      for $k=j+d-1$, where $0\leq j \leq \left\lfloor \frac{d_{1}-1}{2} \right\rfloor$.	      
	\item If $d_{1}$ is even then
	      $|J_{k}| = \frac{3d_{1}^{2}+2d_{1}+4d^{2}+4d-4dd_{1}}{8} + \frac{j(2d+d_{1})}{2} + \frac{3j(j+1)}{2}$,
	      for $k = j + \frac{2d+d_{1}-2}{2}$, where $0\leq j \leq \frac{d_{1}-2}{2}$.
	      
	      If $d_{1}$ is odd then
	      $|J_{k}| = \frac{3d_{1}^{2}+4d^{2}-4dd_{1}-3}{2} + \frac{j(2d+d_{1})}{2} + \frac{3j^{2}}{2}$
	      for $k = j + \frac{2d+d_{1}-3}{2}$, where $0\leq j \leq \frac{d_{1}-1}{2}$.		
  \item $|J_{k}| = \frac{d(d-1)}{2} + d_{1}(d_{1}-1) + j(2d_{1}+d) + \frac{j(j-1)}{2}$, for
        $k=j +d_{1}+d -2$ , where $0  \leq j \leq d-d_{1}$.                 
  \item $|J_{k}| =  \frac{	2d(2d-1) - d_{1}(d_{1}-1)}{2} + j(2d+d_{1}))$, for
        $k=j+2d-2$, where $0 \leq j \leq d_{1}-1$.        
  \item $J_{k} = S_{k}$, for $k\geq 3d-2$.
\end{enumerate}
\end{cor}

\begin{prop}
\em Let $2\leq d_1<d_2=d_3=:d$ be positive integers. Let $f_{1},f_{2},f_{3}\in K[x_{1},x_{2},x_{3}]$ be a regular sequence of homogeneous polynomials of degrees $d_{1},d_{2},d_{3}$. Let $I=(f_{1},f_{2},f_{3})$, $J=\Gin(I)$, the generic initial ideal with respect to the reverse lexicographic order and
$S/I$ has (SLP), then if $d_1$ is even we have:
	\[J=(x_{1}^{d_{1}}, x_{1}^{d_{1}-2j+1}x_{2}^{d-d_{1}-2+3j}, x_{1}^{d_{1}-2j}x_{2}^{d-d_{1}-1+3j}\;for\; 
	     1\leq j \leq  \frac{d_{1}-2}{2}, \]
	\[  x_{1}x_{2}^{\frac{d_{1}+2d-4}{2}}, x_{2}^{\frac{d_{1}+2d-4}{2}},x_{3}^{2j}x_{1}^{2j-1}x_{2}^{\frac{d_{1}+2d}{2} -3j}, \ldots, x_{3}^{2j}x_{2}^{\frac{d_{1}+2d-2}{2}-j}
	    \;for\; 1 \leq j\leq \frac{d_{1}-4}{2},\]
	\[x_{3}^{d_{1}-2+2j}x_{2}^{d-d_{1}-j+1}\{x_{1},x_{2}\}^{d_{1}-1} \;for\;1\leq j \leq d-d_{1},\]	    
  \[  x_{3}^{2d-d_{1}-2+2j}\{x_{1},x_{2}\}^{d_{1}-j}\;for\; 1\leq j\leq d_{1}).\]
	Otherwise, if $d_1$ is odd, then:
	\[J=(x_{1}^{d_{1}}, x_{1}^{d_{1}-2j+1}x_{2}^{d-d_{1}-2+3j}, x_{1}^{d_{1}-2j}x_{2}^{d-d_{1}-1+3j}\;for\; 
	     1\leq j \leq  \frac{d_{1}-1}{2},  \]
	\[ x_{2}^{\frac{d_{1}+2d-1}{2}}, x_{3}x_{2}^{\frac{d_{1}+2d-3}{2}}, x_{3}^{2j+1}x_{1}^{2j}x_{2}^{\frac{d_{1}+2d-3}{2}
	    -3j}, \ldots, x_{3}^{2j+1}x_{2}^{\frac{d_{1}+2d-3}{2}-3j}\;for\; 1 \leq j\leq \frac{d_{1}-3}{2}, \]
	\[x_{3}^{d_{1}-2+2j}x_{2}^{d-d_{1}-j+1}\{x_{1},x_{2}\}^{d_{1}-1} \;for\;1\leq j \leq d-d_{1},\]	    
  \[  x_{3}^{2d-d_{1}-2+2j}\{x_{1},x_{2}\}^{d_{1}-j}\;for\; 1\leq j\leq d_{1}).\]	
\end{prop}

\begin{proof}
We have $|J_{d_1}|=1$, hence $J_{d_1}=\{x_{1}^{d_1}\}$, since $J$ is a strongly stable. Therefore:
\[ Shad(J_{d_1}) = \{x_{1}^{d_1}\{x_{1},x_{2}\} ,x_{1}^{d_1}x_{3}\}. \]
Assume $d>d_1+1$. Using the formulae from $3.12$ we get $|J_{d_{1}+1}|-|Shad(J_{d_{1}})|=0$ so $J_{d_{1}+1} = Shad(J_{d_{1}})$. We show by induction on $1\leq j\leq d-d_{1}$ that 
\[ J_{d_{1}+j} = Shad(J_{d_{1}-1+j}) = \{x_{1}^{d_{1}}\{x_{1},x_{2}\}^{j}, x_{1}^{d_{1}}x_{3}\{x_{1},x_{2}\}^{j-1}, \ldots, x_{1}^{d_{1}}x_{3}^{j}\}.\]
We already prove this for $j=1$. Suppose the assertion is true for some $j<d-d_1$. Since $|J_{d_{1}+j}|=|Shad(J_{d_{1}-1+j})|$ we have $J_{d_{1}+j} = Shad(J_{d_{1}-1+j})$ thus the induction step is done.
In particular, we get:
\[ J_{d-1} = \{x_{1}^{d_{1}}\{x_{1},x_{2}\}^{d-d_{1}-1}, x_{1}^{d_{1}}x_{3}\{x_{1},x_{2}\}^{d-d_{1}-2}, \ldots, x_{1}^{d_{1}}x_{3}^{d-d_{1}-1}\}. \]
which is the same expression as in the case $d=d_1+1$.

In the following, we consider two possibilities. First, suppose $d_1=2$. We have $|J_{d}|-|Shad(J_{d-1})|=2$ so we must add two new generators to $Shad(J_{d-1})$ in order to obtain $J_{d}$. Since $J$ is strongly stable and $x_3$ is a strong Lefschetz element for $S/J$ these new generators
are $x_1 x_2^{d-1}$ and $x_2^{d}$.

Suppose now $d_1>2$. We have $|J_{d}|-|Shad(J_{d-1})|=2$ so we must add
two new generators to $Shad(J_{d-1})$ in order to obtain $J_{d}$. Since $J$ is strongly stable and $x_3$ is a strong Lefschetz element for $S/J$ these new generators
are $x_{1}^{d_{1}-1}x_{2}^{d-d_{1}+1},x_{1}^{d_{1}-2}x_{2}^{d-d_{1}+2}$. Therefore
\[ J_{d} = \{x_{1}^{d_{1}-2}\{x_{1},x_{2}\}^{d-d_{1}}, x_{1}^{d_{1}}x_{3}\{x_{1},x_{2}\}^{d-d_{1}-1}, \ldots, x_{1}^{d_{1}}x_{3}^{d-d_{1}}\}.\]
We prove by induction on $1\leq j \leq \left\lfloor \frac{d_{1}-1}{2} \right\rfloor$ that:
\[ J_{d-1+j} = Shad(J_{d-2+j})\cup \{x_{1}^{d_{1}-2j+1}x_{2}^{d-d_{1}-2+3j}, x_{1}^{d_{1}-2j}x_{2}^{d-d_{1}-1+3j}\} = \]
\[ = \{ x_{1}^{d_{1}-2j}\{x_{1},x_{2}\}^{d-d_{1}-1-3j}, x_{1}^{d_{1}-2j+2}x_{3}\{x_{1},x_{2}\}^{d-d_{1}-4-3j},\ldots,
x_{1}^{d_{1}}x_{3}^{j}\{x_{1},x_{2}\}^{d-d_{1}-1},\]
\[ x_{1}^{d_{1}}x_{3}^{j+1}\{x_{1},x_{2}\}^{d-d_{1}-2},\ldots, x_{1}^{d_{1}}x_{3}^{d-d_{1}+j-1}\}. \]
We already done the case $j=1$. Suppose the assertion is true for some $j<\left\lfloor \frac{d_{1}-1}{2} \right\rfloor$. Since $|J_{d+j}|-|Shad(J_{d-1+j})|=2$ we must add two generators to $Shad(J_{d-1+j})$ in order
to obtain $|J_{d+j}|$ and they are $x_{1}^{d_{1}-2j-1}x_{2}^{d-d_{1}+3j+1}, x_{1}^{d_{1}-2j-2}x_{2}^{d-d_{1}+3j+2}$
because $J$ is strongly stable and $x_3$ is a strong Lefschetz element for $S/J$. Therefore, the induction step is done.

In the following, we consider two possibilities: $d_1$ is even or $d_1$ is odd. First, suppose $d_1$ is even. We have
\[ J_{\frac{d_{1}+2d-4}{2}} = \{x_{1}^{2}\{x_{1},x_{2}\}^{\frac{d_{1}+2d-8}{2}}, x_{1}^{4}x_{3}\{x_{1},x_{2}\}^{\frac{d_{1}+2d-14}{2}}, \ldots, x_{1}^{d_{1}}x_{3}^{\frac{d_{1}-2}{2}} \{x_{1},x_{2}\}^{d-d_{1}-1},\] \[x_{1}^{d_{1}}x_{3}^{\frac{d_{1}}{2}} \{x_{1},x_{2}\}^{d-d_{1}-2}, \ldots,
x_{1}^{d_{1}}x_{3}^{\frac{2d-d_{1}-4}{2}}\}. \]
Since $|J_{\frac{d_{1}+2d-2}{2}}| - |Shad(J_{\frac{d_{1}+2d-4}{2}})| = 2$ we need to add two new monomials to 
$Shad(J_{\frac{d_{1}+2d-4}{2}})$ and since $J$ is strongly stable and $x_3$ is a strong Lefschetz element for $S/J$, they are 
 $x_{1}x_{2}^{\frac{d_{1}+2d-4}{2}}, x_{2}^{\frac{d_{1}+2d-4}{2}}$, thus:
\[ J_{\frac{d_{1}+2d-2}{2}} =\{\{x_{1},x_{2}\}^{\frac{d_{1}+2d-2}{2}},x_{3}x_{1}^{2}\{x_{1},x_{2}\}^{\frac{d_{1}+2d-8}{2}}, x_{1}^{4}x_{3}^{2}\{x_{1},x_{2}\}^{\frac{d_{1}+2d-14}{2}}, \ldots, \] \[x_{1}^{d_{1}}x_{3}^{\frac{d_{1}}{2}} \{x_{1},x_{2}\}^{d-d_{1}-1},x_{1}^{d_{1}}x_{3}^{\frac{d_{1}+2}{2}} \{x_{1},x_{2}\}^{d-d_{1}-2}, \ldots,
x_{1}^{d_{1}}x_{3}^{\frac{2d-d_{1}-2}{2}}\}. \]

One can easily show by induction on $1 \leq j\leq \frac{d_{1}-2}{2}$ that
\[ J_{\frac{d_{1}+2d-2}{2} + j} = Shad(J_{\frac{d_{1}+2d-4}{2} + j}) \cup \{x_{3}^{2j}x_{1}^{2j-1}x_{2}^{\frac{d_{1}+2d}{2} -3j}, \ldots, x_{3}^{2j}x_{2}^{\frac{d_{1}+2d-2}{2}-j}\} = \]
\[ = \{ \{x_{1},x_{2}\}^{\frac{d_{1}+2d-2}{2} + j}, x_{3}\{x_{1},x_{2}\}^{\frac{d_{1}+2d-4}{2} + j}, \ldots, x_{3}^{2j}\{x_{1},x_{2}\}^{\frac{d_{1}+2d-2}{2} - j}, \]\[
x_{3}^{2j+1}x_{1}^{2j+2}\{x_{1},x_{2}\}^{\frac{d_{1}+2d-8}{2} - 3j}, \ldots, x_{3}^{\frac{d_{1}}{2}}x_{1}^{d_{1}}\{x_{1},x_{2}\}^{d-d_{1}-1},\ldots,x_{1}^{d_{1}}x_{3}^{\frac{2d-d_{1}-2}{2}+j}\}.\]
The assertion was already done for $j=1$ and the induction step is similar.

If $d_1$ is odd, we get
\[ J_{\frac{d_{1}+2d-3}{2}} = \{x_{1}\{x_{1},x_{2}\}^{\frac{d_{1}+2d-5}{2}}, x_{1}^{3}x_{3}\{x_{1},x_{2}\}^{\frac{d_{1}+2d-11}{2}}, \ldots, x_{1}^{d_{1}}x_{3}^{\frac{d_{1}-1}{2}} \{x_{1},x_{2}\}^{d-d_{1}-1},\] \[x_{1}^{d_{1}}x_{3}^{\frac{d_{1}+1}{2}} \{x_{1},x_{2}\}^{d-d_{1}-2}, \ldots,
x_{1}^{d_{1}}x_{3}^{\frac{2d-d_{1}-3}{2}}\}.\]

Since $|J_{\frac{d_{1}+2d-1}{2}}| - |Shad(J_{\frac{d_{1}+2d-3}{2}})|  = 2$ we add two generators to 
$Shad(J_{\frac{d_{1}+2d-3}{2}})$ in order to obtain $J_{\frac{d_{1}+2d-1}{2}}$ and they must be
$x_{2}^{\frac{d_{1}+2d-1}{2}}, x_{3}x_{2}^{\frac{d_{1}+2d-3}{2}}$, therefore:
\[ J_{\frac{d_{1}+2d-1}{2}} =\{\{x_{1},x_{2}\}^{\frac{d_{1}+2d-1}{2}}, x_{3}\{x_{1},x_{2}\}^{\frac{d_{1}+2d-3}{2}}, x_{1}^{3}x_{3}^{2}\{x_{1},x_{2}\}^{\frac{d_{1}+2d-11}{2}}, \ldots, \]\[ x_{1}^{d_{1}}x_{3}^{\frac{d_{1}+1}{2}} \{x_{1},x_{2}\}^{d-d_{1}-1},x_{1}^{d_{1}}x_{3}^{\frac{d_{1}+3}{2}} \{x_{1},x_{2}\}^{d-d_{1}-2}, \ldots,
x_{1}^{d_{1}}x_{3}^{\frac{2d-d_{1}-1}{2}}\}.\]
Since $|J_{\frac{d_{1}+2d+1}{2}}| - |Shad(J_{\frac{d_{1}+2d-1}{2}})|  = 3$, we add $3$ new generators to 
$Shad(J_{\frac{d_{1}+2d-1}{2}})$ in order to obtain $J_{\frac{d_{1}+2d+1}{2}}$ and they are
 $x_{1}^{2}x_{3}^{3}x_{2}^{\frac{d_{1}+2d-9}{2}}, x_{1}x_{3}^{3}x_{2}^{\frac{d_{1}+2d-7}{2}}, x_{3}^{3}x_{2}^{\frac{d_{1}+2d-5}{2}}$, therefore
\[ J_{\frac{d_{1}+2d+1}{2}} = \{\{x_{1},x_{2}\}^{\frac{d_{1}+2d+1}{2}}, x_{3}\{x_{1},x_{2}\}^{\frac{d_{1}+2d-1}{2}}, x_{3}^{2}\{x_{1},x_{2}\}^{\frac{d_{1}+2d-3}{2}},x_{3}^{3}\{x_{1},x_{2}\}^{\frac{d_{1}+2d-5}{2}}, \ldots,\] 
\[ x_{1}^{d_{1}}x_{3}^{\frac{d_{1}+3}{2}} \{x_{1},x_{2}\}^{d-d_{1}-1},x_{1}^{d_{1}}x_{3}^{\frac{d_{1}+5}{2}} \{x_{1},x_{2}\}^{d-d_{1}-2}, \ldots, x_{1}^{d_{1}}x_{3}^{\frac{2d-d_{1}+1}{2}}\}. \]

One can easily prove by induction on $1\leq j\leq \frac{d_{1}-3}{2}$ that
\[ J_{\frac{d_{1}+2d-1}{2} + j} = Shad(J_{\frac{d_{1}+2d-3}{2} + j})\cup \{ x_{3}^{2j+1}x_{1}^{2j}x_{2}^{\frac{d_{1}+2d-3}{2} -3j}, \ldots, x_{3}^{2j+1}x_{2}^{\frac{d_{1}+2d-3}{2}-3j}\} = \]
\[ = \{ \{x_{1},x_{2}\}^{\frac{d_{1}+2d-1}{2} + j}, x_{3}\{x_{1},x_{2}\}^{\frac{d_{1}+2d-3}{2} + j}, \ldots, x_{3}^{2j+2}\{x_{1},x_{2}\}^{\frac{d_{1}+2d-3}{2} - j},\]
\[ x_{3}^{2j+2}x_{1}^{2j+3}\{x_{1},x_{2}\}^{\frac{d_{1}+2d-11}{2} - 3j}, \ldots, x_{3}^{\frac{d_{1}+1}{2}}x_{1}^{d_{1}}\{x_{1},x_{2}\}^{d-d_{1}-1},\ldots,x_{1}^{d_{1}}x_{3}^{\frac{2d-d_{1}-1}{2}+j}\}.\]
The assertion was already proved for $j=1$ and the induction step is similar.

In all cases, we obtain
\[ J_{d_{1}+d-2} = \{ \{x_{1},x_{2}\}^{d_{1}+d-2}, x_{3}\{x_{1},x_{2}\}^{d_{1}+d-3}, \ldots, x_{3}^{d_{1}-2}\{x_{1},x_{2}\}^{d},\] \[ x_{1}^{d_{1}}x_{3}^{d_{1}-1}\{x_{1},x_{2}\}^{d-d_{1}-1},\ldots, x_{1}^{d_{1}}x_{3}^{d-2}\}. \]
We have $|J_{d_{1}+d-1}|-|Shad(J_{d_{1}+d-2})|=d_1$ so we must add $d_1$ new monomials to $Shad(J_{d_{1}+d-2})$
in order to obtain $J_{d_{1}+d-1}$. Since $J$ is strongly stable and $x_3$ is a strong Lefschetz element for $S/J$, these new monomials are 
 $x_{3}^{d_{1}}x_{2}^{d-d_{1}}\{x_{1},x_{2}\}^{d_{1}-1}$, therefore
\[ J_{d_{1}+d-1} = ( \{x_{1},x_{2}\}^{d_{1}+d-1}, x_{3}\{x_{1},x_{2}\}^{d_{1}+d-2}, \ldots, x_{3}^{d_{1}}\{x_{1},x_{2}\}^{d-1},\] \[ x_{1}^{d_{1}}x_{3}^{d_{1}+1}\{x_{1},x_{2}\}^{d-d_{1}-2},\ldots, x_{1}^{d_{1}}x_{3}^{d-1}).\]

One can easily prove by induction on $1\leq j \leq d-d_{1}$ that
\[ J_{d_{1}+d-1+j} = Shad(J_{d_{1}+d-2+j}) \cup \{ x_{3}^{d_{1}-2+2j}x_{2}^{d-d_{1}-j+1}\{x_{1},x_{2}\}^{d_{1}-1}\} = \]
\[ = \{ \{x_{1},x_{2}\}^{d_{1}+d-2+j}, x_{3}\{x_{1},x_{2}\}^{d_{1}+d-3+j}, \ldots, x_{3}^{d_{1}-2+2j}\{x_{1},x_{2}\}^{d-j},\] \[ x_{1}^{d_{1}}x_{3}^{d_{1}-1+2j}\{x_{1},x_{2}\}^{d-d_{1}-j-1},\ldots, x_{1}^{d_{1}}x_{3}^{d-2+j}\}  \]
the case $j=1$ being already done and than, the induction step being similar. In particular,
$J_{2d-2} = \{\{x_{1},x_{2}\}^{2d-2},x_{3}\{x_{1},x_{2}\}^{2d-3},\ldots ,x_{3}^{2d-d_{1}-2}\{x_{1},x_{2}\}^{d_{1}}\}$.

We have $|J_{2d-1}| -|Shad(J_{2d-2})|=d_1$ so we must add $d_1$ new generators to $Shad(J_{2d-2})$ in order to
obtain $J_{2d-1}$. Since $J$ is strongly stable and $x_3$ is a strong Lefschetz element for $S/J$, these new monomials are $x_{3}^{2d-d_{1}}\{x_{1},x_{2}\}^{d_{1}-1}$ so
\[J_{2d-1} = \{\{x_{1},x_{2}\}^{2d-1},x_{3}\{x_{1},x_{2}\}^{2d-2},\ldots ,x_{3}^{2d-d_{1}}\{x_{1},x_{2}\}^{d_{1}-1}\}. \]
One can easily show by induction on $1\leq j\leq d_{1}$ that
\[J_{2d-2+j} = Shad(J_{2d-3+j})\cup \{x_{3}^{2d-d_{1}-2+2j}\{x_{1},x_{2}\}^{d_{1}-j}\} = \]
 \[ = \{\{x_{1},x_{2}\}^{2d-2+j},\ldots,x_{3}^{2d-d_{1}-2+2j}\{x_{1},x_{2}\}^{d_{1}-j}\}.\]
Finally, we obtain $J_{d_{1}+2d-2}=S_{d_{1}+2d-2}$ and therefore we cannot add new minimal generators of $J$ in degrees $>d_{1}+2d-2$.
\end{proof}

\begin{cor}
In the conditions of the above proposition, the number of minimal generators of $J$ is $d_1(d+1) - 
\left(\frac{d_1}{2} \right)^{2} + 1$ if $d$ is even; $d_1(d+1) - \frac{d_1^{2}-1}{4} + 1$ if $d$ is odd.
\end{cor}

\begin{ex}
If $d_{1}=4$, $d_{2}=d_{3}=6$, then $J = (x_1^{4},\; x_1^{3}x_2^{3},\; x_1^{2}x_2^{4},\;x_2^{7},\; x_3x_2^{6},\;x_3^{2}x_{1}x_{2}^{5},$ \linebreak 
$ x_3^{2}x_{2}^{6},\; x_3^{4}x_{2}^{2}\{x_1,x_2\}^{3},\;  
	     x_3^{6}x_{2}\{x_1,x_2\}^{3},\; x_3^{8}\{x_1,x_2\}^{3},\; x_3^{10}\{x_1,x_2\}^{2},\; 
	       x_3^{12}\{x_1,x_2\},\; x_3^{14})$. \vspace{5 pt}

If $d_{1}=3$ and $d_{2}=d_{3}=6$, then: $J = (x_1^{3},\; x_1^{2}x_{2}^{4},\; x_1x_2^{5},\; x_1x_2^{6},\;x_2^{7},\;x_3^{3}x_2^{3}\{x_1,x_2\}^{2},\;$
\linebreak $ x_3^{5}x_2^{2}\{x_1,x_2\}^{2},\;   
	       x_3^{7}x_2\{x_1,x_2\}^{2},\;,x_3^{9}\{x_1,x_2\}^{2},\; x_3^{11}\{x_1,x_2\}^,\; x_3^{11}  )$.       
\end{ex}

\medskip
\begin{itemize}
\item Subcase $d_1<d_2<d_3$.
\end{itemize} 
\medskip

\begin{prop}
\em Let $2\leq d_1<d_2<d_3$ be positive integers such that \linebreak $d_{1}+d_{2}>d_{3}+1$. The Hilbert function of the standard graded complete intersection $A=K[x_{1},x_{2},x_{3}]/I$, where $I$ is the ideal generated by $f_{1}$, $f_{2}$, $f_{3}$, with $f_i$ homogeneous polynomials of degree $d_i$, for all $i$, with $1\leq i\leq 3$, has the form:  
\begin{enumerate}
	\item $H(A,k) = \binom{k+2}{k}$, for $k\leq d_{1}-2$.
	\item $H(A,k) = \binom{d_{1}+1}{2} + jd_{1}$, for $k=j+d_{1}-1$, where $0\leq j \leq d_{2}-d_{1}$.		
	\item $H(A,k) = \binom{d_{1}+1}{2} + d_{1}(d_{2}-d_{1}) + \sum_{i=1}^{j}(d_{1}-i)$, 
	      for $k=j+d_{2}-1$, where $0 \leq k \leq d_{3}-d_{2}$.	      
	\item $H(A,k) = \binom{d_{1}+1}{2} + d_{1}(d_{2}-d_{1}) + \sum_{i=1}^{d_{3}-1}(d_{1}-j) + \sum_{i=1}^{j}
	               (d_{1}+d_{2}-d_{3}-2i)$, for $k=j+d_{3}-1$, where $0\leq j \leq [\frac{d_{1}+d_{2}-d_{3}-1}{2}]$.
	\item $H(A,k) = H(A,d_{1}+d_{2}+d_{3}-3-k)$ for $k>d_{3}-1 + [\frac{d_{1}+d_{2}-d_{3}-1}{2}]$.
\end{enumerate}
\end{prop}

\begin{proof}
It follows from \cite[Lemma 2.9(b)]{HPV}.
\end{proof}

\begin{cor}
\em In the conditions of Proposition $3.16$, let $J=\Gin(I)$ be the generic initial ideal of $I$ with respect to the reverse lexicographic order. If we denote by $J_{k}$ the set of monomials of $J$ of degree $k$, then:
\begin{enumerate}	
	\item $|J_{k}| = 0$, for $k\leq d_{1}-1$.
	\item $|J_{k}| = j(j+1)/2$, for $k = j +d_{1} - 1$, where $0 \leq j \leq d_{2}-d_{1}$.	
	\item $|J_{k}| = d_{2}(d_{2}-1) + j(d_{2}-d_{1}) + j(j+1)$,
	      for $k = j +d_{2} - 1$, where $0 \leq j \leq d_{3}-d_{2}$.	      
	\item $|J_{k}| = |J_{d_{3}-1}|+j(2d_{3}-d_{1}-d_{2})+\frac{3j(j+1)}{2}$, for $k = j+d_{3}-1$, 
	      where $0 \leq j \leq [\frac{\alpha-1}{2}]$.	      
	\item If $d_{1}+d_{2}+d_{3}$ is even then 
	      $|J_{k}| = |J_{\frac{d_{1}+d_{2}+d_{3}-4}{2}}|  + \frac{(j+1)(d_{1}+d_{2}+d_{3})}{2} + \frac{3j(j+1)}{2}$,
	      for $k = j + \frac{d_{1}+d_{2}+d_{3}-2}{2}$, where $0 \leq j \leq \frac{d_{1}+d_{2}-d_{3}-2}{2}$.
	
	      If $d_{1}+d_{2}+d_{3}$ is odd then 
	      $|J_{k}| = |J_{\frac{d_{1}+d_{2}+d_{3}-3}{2}}| +  \frac{j(d_{1}+d_{2}+d_{3})}{2} + \frac{3j^{2}}{2}$, 
	      for $k = j + \frac{d_{1}+d_{2}+d_{3}-3}{2}$, where $0 \leq j \leq \frac{d_{1}+d_{2}-d_{3}-1}{2}$.
  \item $|J_{k}| = |J_{d_{1}+d_{2}-2}| + j(2d_{1}+2d_{2}-d_{3}-1)+j^{2}$, 
        for $k=j+d_{1}+d_{2}-2$, where $0 \leq d_{3}-d_{2}$.       
  \item $|J_{k}| = |J_{d_{1}+d_{3}-2}|+j(2d_{1}+d_{3}-1)+\frac{j(j+1)}{2}$, for $k = j+d_{1}+d_{3}-2$,
        where $0 \leq j \leq d_{2}-d_{1}$.        
  \item $|J_{k}| = |J_{d_{2}+d_{3}-1}|+j(d_{1}+d_{2}+d_{3})$, for $k = j + d_{1}+d_{3}-2$,
        where $0 \leq j \leq d_{1}-1$        
  \item $J_{k} = S_{k}$, for $k\geq 3d-2$.
\end{enumerate}
\end{cor}

\pagebreak

\begin{prop}
\em Let $2\leq d_1<d_2<d_3$ be positive integers such that $d_{1}+d_{2} > d_3+1$. Let $f_{1},f_{2},f_{3}\in K[x_{1},x_{2},x_{3}]$ be a regular sequence of homogeneous polynomials of degrees $d_{1},d_{2},d_{3}$. Let
$\alpha =d_{1}+d_{2}-d_{3}$. Let $I=(f_{1},f_{2},f_{3})$, $J=\Gin(I)$, the generic initial ideal with respect to the reverse lexicographic order, and suppose $S/I$ has (SLP). If $\alpha$ is even, then:
\[J = (x_{1}^{d_{1}},x_{1}^{d_{1}-1}x_{2}^{d_{2}-d_{1}+1},x_{1}^{d_{1}-2}x_{2}^{d_{2}-d_{1}+3},\ldots, 
        x_{1}^{d_{1}+d_{2}-d_{3}}x_{2}^{2d_{3}-d_{1}-d_{2}-1}, \] 
\[  x_{1}^{d_{1}+d_{2}-d_{3}-2j}x_{2}^{2d_{3}-d_{1}-d_{2}+3j-1},
x_{1}^{d_{1}+d_{2}-d_{3}-2j+1}x_{2}^{2d_{3}-d_{1}-d_{2}+3j-2} 
 for\;j=1,\ldots,\frac{\alpha-2}{2},\]
\[x_{2}^{\frac{d_{1}+d_{2}+d_{3}-2}{2}},x_{1}x_{2}^{\frac{d_{1}+d_{2}+d_{3}-4}{2}}, 
x_{3}^{2j}x_{2}^{\frac{d_{1}+d_{2}+d_{3}}{2}-3j} \{x_{1},x_{2}\}^{2j-1} for\;j=1,\ldots,\frac{\alpha-2}{2}, \]
\[ x_{1}^{d_{1}+d_{2}-d_{3}+j-2}x_{2}^{2d_{3}-d_{1}-d_{2}-2j+2}x_{3}^{d_{1}+d_{2}-d_{3}+2j-2},..,
  x_{2}^{d_{3}-j}x_{3}^{d_{1}+d_{2}-d_{3}+2j-2}\; for\;j=1,..,d_{3}-d_{2},  \]
\[  x_{1}^{d_{1}-1}x_{3}^{d_{1}+d_{3}-d_{2}+2j-2}x_{2}^{d_{2}-d_{1}-j+1} , \ldots, x_{3}^{d_{1}+d_{3}-d_{2}+2j-2}x_{2}^{d_{2}-j}\; for\; j=1,\ldots,d_{2}-d_{1}, \]
\[ \{x_{1},x_{2}\}^{d_{1}-j}x_{3}^{d_{2}+d_{3}-d_{1}-2+2j}\; for\; j=1,\ldots,d_{1}). \]

Otherwise, if $\alpha$ is odd, then:
\[J = (x_{1}^{d_{1}},x_{1}^{d_{1}-1}x_{2}^{d_{2}-d_{1}+1},x_{1}^{d_{1}-2}x_{2}^{d_{2}-d_{1}+3},\ldots, 
        x_{1}^{d_{1}+d_{2}-d_{3}}x_{2}^{2d_{3}-d_{1}-d_{2}-1},\]
\[ x_{1}^{d_{1}+d_{2}-d_{3}-2j}x_{2}^{2d_{3}-d_{1}-d_{2}+3j-1},
x_{1}^{d_{1}+d_{2}-d_{3}-2j+1}x_{2}^{2d_{3}-d_{1}-d_{2}+3j-2} ,
 j=1,\ldots,\frac{\alpha-1}{2}, \]
\[x_{2}^{\frac{d_{1}+d_{2}+d_{3}-1}{2}},x_{3}x_{2}^{\frac{d_{1}+d_{2}+d_{3}-3}{2}} ,x_{1}^{2j}x_{3}^{2j+1}x_{2}^{\frac{d_{1}+d_{2}+d_{3}-3}{2} -3j},.., x_{3}^{2j+1}x_{2}^{\frac{d_1+d_2+d_3-3}{2} - j} ,j=1,..,\frac{\alpha-3}{2}, \]
\[ x_{1}^{d_{1}+d_{2}-d_{3}+j-2}x_{2}^{2d_{3}-d_{1}-d_{2}-2j+2}x_{3}^{d_{1}+d_{2}-d_{3}+2j-2},..,
  x_{2}^{d_{3}-j}x_{3}^{d_{1}+d_{2}-d_{3}+2j-2}\; for\;j=1,..,d_{3}-d_{2},\]
\[   x_{1}^{d_{1}-1}x_{3}^{d_{1}+d_{3}-d_{2}+2j-2}x_{2}^{d_{2}-d_{1}-j+1} , \ldots, x_{3}^{d_{1}+d_{3}-d_{2}+2j-2}x_{2}^{d_{2}-j}\; for\; j=1,\ldots,d_{2}-d_{1}, \]
\[ \{x_{1},x_{2}\}^{d_{1}-j}x_{3}^{d_{2}+d_{3}-d_{1}-2+2j}\; for\; j=1,\ldots,d_{1}). \]
\end{prop}

\begin{proof}
We have $|J_{d_1}|=1$, hence $J_{d_1}=\{x_{1}^{d_1}\}$,
since $J$ is strongly stable. Therefore:
\[ Shad(J_{d_1}) = \{x_{1}^{d_1}\{x_{1},x_{2}\} ,x_{1}^{d_1}x_{3}\}. \]
Assume $d_2>d_1+1$. Using the formulae from $3.17$ we get $|J_{d_1+1}|-|Shad(J_{d_1})|=0$, therefore
\[J_{d_1+1} = Shad(J_{d_1}) = \{x_{1}^{d_1}\{x_{1},x_{2}\} ,x_{1}^{d_1}x_{3}\}. \]
Using induction on $1\leq j\leq d_2-d_1-1$ we prove that 
\[ J_{d_{1}+j} = Shad(J_{d_1+j-1}) = \{ x_{1}^{d_{1}}\{x_{1},x_{2}\}^{j} ,\ldots,x_{3}^{j}x_{1}^{d_{1}}\}. \]
Indeed, this assertion was already proved for $j=1$, and if we suppose that is true for some $j<d_2-d_1-1$ we
get $|J_{d_{1}+j+1}|=|Shad(J_{d_1+j})|$ so we are done. In particular, we obtain
$ J_{d_2-1} = \{ x_{1}^{d_{1}}\{x_{1},x_{2}\}^{d_{2}-d_{1}-1}, x_{3}x_{1}^{d_{1}}\{x_{1},x_{2}\}^{d_{2}-d_{1}-2},\ldots,x_{3}^{d_{2}-d_{1}-1}x_{1}^{d_{1}} \}.$

We have $|J_{d_2}|-|Shad(J_{d_2-1})|=1$ so we must add a new generator to $Shad(J_{d_2-1})$ in order to obtain
$J_{d_2}$. But since $J$ is strongly stable and $x_3$ is a strong Lefschetz element for $S/J$, this new generator is $x_{1}^{d_{1}-1}x_{2}^{d_{2}-d_{1}+1}$ and therefore 
$ J_{d_{2}} = \{ x_{1}^{d_{1}-1}\{x_{1},x_{2}\}^{d_{2}-d_{1}+1}, x_{3}x_{1}^{d_{1}}\{x_{1},x_{2}\}^{d_{2}-d_{1}-1},\ldots,x_{3}^{d_{2}-d_{1}}x_{1}^{d_{1}}\}. $

We show by induction on $1\leq j\leq d_3-d_2$ that 
\[ J_{d_2-1+j} = Shad(J_{d_2-2+j})\cup \{x_{1}^{d_{1}-j} x_{2}^{d_{2}-d_{1}+2j-1}\} = \]
\[= \{ x_{1}^{d_{1}-j}\{x_{1},x_{2}\}^{d_{2}-d_{1}+2j-1}, \ldots, x_{3}^{j}x_{1}^{d_{1}}\{x_{1},x_{2}\}^{d_{2}-d_{1}-1},\ldots,x_{3}^{d_{2}-d_{1}+j-1}x_{1}^{d_{1}}\}.\]
The first step of induction was already done. Suppose the assertion is true for some $j<d_{3}-d_{2}$. Since
$|J_{d_2+j}|-|Shad(J_{d_2-1+j})|=1$ and $J$ is strongly stable and $x_3$ is a strong Lefschetz element for $S/J$, we add $x_{1}^{d_{1}-j-1}x_{2}^{d_{2}-d_{1}+2j+1}$ 
to $Shad(J_{d_2-1+j})$ in order to obtain $J_{d_2+j}$. Thus, we are done. In particular, we get
\[ J_{d_{3}-1} = \{x_{1}^{d_{1}+d_{2}-d_{3}}\{x_{1},x_{2}\}^{2d_{3}-d_{1}-d_{2}-1}, x_{3}x_{1}^{d_{1}+d_{2}-d_{3}+1}\{x_{1},x_{2}\}^{2d_{3}-d_{1}-d_{2}-3}, \ldots,\] \[ x_{3}^{d_{3}-d_{2}}x_{1}^{d_{1}}\{x_{1},x_{2}\}^{d_{2}-d_{1}-1},\ldots,x_{3}^{d_{3}-d_{1}-1}x_{1}^{d_{1}}\}.\]

We consider first $\alpha=2$, i.e. $d_3 = d_1 + d_2 - 2$. In this case, since 
$|J_{d_3}|-|Shad(J_{d_3 -1})|=2$ we add two new generators to $Shad(J_{d_3-1})$ in order to obtain $J_{d_3}$.
Since $J$ is strongly stable and $x_3$ is a strong Lefschetz element for $S/J$, these new generators are $x_{2}^{d_{3}}$ and $x_{1}x_{2}^{d_{3}-1}$.

Suppose now $\alpha\geq 2$. We have $|J_{d_3}|-|Shad(J_{d_3 -1})|=2$ so we must add two new generators
to $Shad(J_{d_3-1})$ to obtain $J_{d_3}$. Since $J$ is strongly stable and $x_3$ is a strong Lefschetz element for $S/J$, these new generators are $x_{1}^{d_{1}+d_{2}-d_{3}-2}x_{2}^{2d_{3}-d_{1}-d_{2}+2}$ and
$x_{1}^{d_{1}+d_{2}-d_{3}-1}x_{2}^{2d_{3}-d_{1}-d_{2}+1}$, therefore
\[J_{d_{3}}= \{x_{1}^{d_{1}+d_{2}-d_{3}-2}\{x_{1},x_{2}\}^{2d_{3}-d_{1}-d_{2}+2}, x_{3}x_{1}^{d_{1}+d_{2}-d_{3}}\{x_{1},x_{2}\}^{2d_{3}-d_{1}-d_{2}-1},\]\[ x_{3}^{2}x_{1}^{d_{1}+d_{2}-d_{3}+1}\{x_{1},x_{2}\}^{2d_{3}-d_{1}-d_{2}-3},\ldots,
x_{3}^{d_{3}-d_{2}}x_{1}^{d_{1}}\{x_{1},x_{2}\}^{d_{2}-d_{1}},\ldots,x_{3}^{d_{3}-d_{1}}x_{1}^{d_{1}}\}.\]
One can prove by induction on $1\leq j\leq \left\lfloor \frac{\alpha-1}{2} \right\rfloor$ that $J_{d_{3}-1+j}$ is the set
\[ Shad(J_{d_{3}-2+j})\cup \{ x_{1}^{d_{1}+d_{2}-d_{3}-2j}x_{2}^{2d_{3}-d_{1}-d_{2}+3j-1}, x_{1}^{d_{1}+d_{2}-d_{3}-2j+1}x_{2}^{2d_{3}-d_{1}-d_{2}+3j-2} \} = \]
\[ \{x_{1}^{d_{1}+d_{2}-d_{3}-2j}\{x_{1},x_{2}\}^{2d_{3}-d_{1}-d_{2}+3j-1},
\ldots, x_{1}^{d_{1}+d_{2}-d_{3}-2}x_{3}^{j-1}\{x_{1},x_{2}\}^{2d_{3}-d_{1}-d_{2}+2},
x_{3}^{j}J_{d_{3}-1}\}.\]
Indeed, the assertion was proved for $j=1$ and the induction step is similar. In the following, we must consider
two possibilities: $\alpha$ is even or $\alpha$ is odd. Suppose first $\alpha$ is even. We obtain that $J_{\frac{d_{1}+d_{2}+d_{3}-4}{2}}$ is the set
\[ \{ x_{1}^{2}\{x_{1},x_{2}\}^{\frac{d_{1}+d_{2}+d_{3}-8}{2}}, \ldots,
x_{1}^{d_{1}+d_{2}-d_{3}-2}x_{3}^{\frac{d_{1}+d_{2}-d_{3}-4}{2}}
\{x_{1},x_{2}\}^{2d_{3}-d_{1}-d_{2}+2}, x_{3}^{\frac{d_{1}+d_{2}-d_{3}-2}{2}}J_{d_{3}-1}\}.\] 
We have $|J_{\frac{d_{1}+d_{2}+d_{3}-2}{2}}|-|Shad(J_{\frac{d_{1}+d_{2}+d_{3}-4}{2}})|=2$ so we must add two
generators to $Shad(J_{\frac{d_{1}+d_{2}+d_{3}-4}{2}})$ and, since $J$ is strongly stable and $x_3$ is a strong Lefschetz element for $S/J$, these new generators are
$x_{2}^{\frac{d_{1}+d_{2}+d_{3}-2}{2}}$ and $x_{1}x_{2}^{\frac{d_{1}+d_{2}+d_{3}-4}{2}}$. Therefore
\[ J_{\frac{d_{1}+d_{2}+d_{3}-2}{2}} = \{ \{x_{1},x_{2}\}^{\frac{d_{1}+d_{2}+d_{3}-2}{2}},
x_{1}^{2}x_{3}\{x_{1},x_{2}\}^{\frac{d_{1}+d_{2}+d_{3}-8}{2}},\ldots,\] 
\[  x_{1}^{d_{1}+d_{2}-d_{3}-2}x_{3}^{\frac{d_{1}+d_{2}-d_{3}-2}{2}}
\{x_{1},x_{2}\}^{2d_{3}-d_{1}-d_{2}+2},x_{3}^{\frac{d_{1}+d_{2}-d_{3}}{2}}J_{d_{3}-1}\}.\]

One can easily prove by induction on $1\leq j \leq \frac{\alpha-4}{2}$ that $J_{\frac{d_{1}+d_{2}+d_{3}-2}{2} + j}$
is the set
\[ = Shad(J_{\frac{d_{1}+d_{2}+d_{3}-4}{2} + j})\cup 
\{ x_{3}^{2j}x_{1}^{2j-1} x_{2}^{\frac{d_{1}+d_{2}+d_{3}-2}{2}-3j},\ldots,x_{3}^{2j}x_{2}^{\frac{d_{1}+d_{2}+d_{3}-2}{2}-j} \} = \]
\[\{ \{x_{1},x_{2}\}^{\frac{d_{1}+d_{2}+d_{3}-2}{2}+j},\ldots,x_{3}^{2j} \{x_{1},x_{2}\}^{\frac{d_{1}+d_{2}+d_{3}-2}{2}-j},x_{1}^{2j+2}x_{3}^{2j+1}\{x_{1},x_{2}\}^{\frac{d_{1}+d_{2}+d_{3}-8}{2}-3j},\ldots,\]
\[ x_{1}^{d_{1}+d_{2}-d_{3}-2}x_{3}^{\frac{d_{1}+d_{2}-d_{3}-2}{2}+j}\{x_{1},x_{2}\}^{2d_{3}-d_{1}-d_{2}+2},
  x_{3}^{\frac{d_{1}+d_{2}-d_{3}}{2}+j}J_{d_{3}-1} \}.\]

Suppose now $\alpha$ is odd. We have that $J_{\frac{d_{1}+d_{2}+d_{3}-3}{2}}$ is the set
\[\{ x_{1}\{x_{1},x_{2}\}^{\frac{d_{1}+d_{2}+d_{3}-5}{2}}, \ldots,
x_{1}^{d_{1}+d_{2}-d_{3}-2}x_{3}^{\frac{d_{1}+d_{2}-d_{3}-3}{2}}\{x_{1},x_{2}\}^{2d_{3}-d_{1}-d_{2}+2},x_{3}^{\frac{d_{1}+d_{2}-d_{3}-1}{2}}J_{d_{3}-1}\}.\]
We have $|J_{\frac{d_{1}+d_{2}+d_{3}-1}{2}}|-| Shad(J_{\frac{d_{1}+d_{2}+d_{3}-3}{2}})|=2$ so we must add two
new generators to $Shad(J_{\frac{d_{1}+d_{2}+d_{3}-3}{2}})$ in order to obtain $J_{\frac{d_{1}+d_{2}+d_{3}-1}{2}}$.
Since $J$ is strongly stable and $x_3$ is a strong Lefschetz element for $S/J$, these new generators are $x_{2}^{\frac{d_{1}+d_{2}+d_{3}-1}{2}}$ and $x_{3}x_{2}^{\frac{d_{1}+d_{2}+d_{3}-3}{2}}$. Therefore
\[J_{\frac{d_{1}+d_{2}+d_{3}-1}{2}} =  \{ \{x_{1},x_{2}\}^{\frac{d_{1}+d_{2}+d_{3}-1}{2}},
x_{3}\{x_{1},x_{2}\}^{\frac{d_{1}+d_{2}+d_{3}-3}{2}}, x_{1}^{3}x_{3}^{2}\{x_{1},x_{2}\}^{\frac{d_{1}+d_{2}+d_{3}-11}{2}}, \ldots, \] 
\[ x_{1}^{d_{1}+d_{2}-d_{3}-2}x_{3}^{\frac{d_{1}+d_{2}-d_{3}-1}{2}}\{x_{1},x_{2}\}^{2d_{3}-d_{1}-d_{2}+2},
x_{3}^{\frac{d_{1}+d_{2}-d_{3}+1}{2}}J_{d_{3}-1} \}.\]

One can easily prove by induction on $1\leq j\leq \frac{\alpha-3}{2}$ that
\[ J_{\frac{d_{1}+d_{2}+d_{3}-1}{2} + j} = Shad(J_{\frac{d_{1}+d_{2}+d_{3}-3}{2} + j})\cup\{ 
x_{1}^{2j}x_{3}^{2j+1}x_{2}^{\frac{d_{1}+d_{2}+d_{3}-3}{2} -3j},\ldots, x_{3}^{2j+1}x_{2}^{\frac{d_{1}+d_{2}+d_{3}-3}{2} - j}\}. \]
For $j=1$, we notice that $|J_{\frac{d_{1}+d_{2}+d_{3}+1}{2}}| - |Shad(J_{\frac{d_{1}+d_{2}+d_{3}-1}{2}})|=3$ so we
must add $3$ new monomials to $Shad(J_{\frac{d_{1}+d_{2}+d_{3}-1}{2}})$ in order to obtain $J_{\frac{d_{1}+d_{2}+d_{3}+1}{2}}$. But, since $J$ is strongly stable and $x_3$ is a strong Lefschetz element for $S/J$, they are exactly 
$x_{1}^{2}x_{3}^{3}x_{2}^{\frac{d_{1}+d_{2}+d_{3}-3}{2} -3}$, $x_{1}x_{3}^{3}x_{2}^{\frac{d_{1}+d_{2}+d_{3}-3}{2}-2}$ and $x_{3}^{3}x_{2}^{\frac{d_{1}+d_{2}+d_{3}-3}{2}-1}$, as required. The induction step is similar.

In all cases, we obtain that $J_{d_{1}+d_{2}-2}$ is the set
\[ \{ \{x_{1},x_{2}\}^{d_{1}+d_{2}-2},x_{3}\{x_{1},x_{2}\}^{d_{1}+d_{2}-3},\ldots, x_{3}^{d_{1}+d_{2}-d_{3}-2}\{x_{1},x_{2}\}^{d_{3}},x_{3}^{d_{1}+d_{2}-d_{3}-1}J_{d_{3}-1}\}.\]
We have $|J_{d_{1}+d_{2}-1}|-|Shad(J_{d_{1}+d_{2}-2})| = \alpha$, so we must add $\alpha$ new generators
to obtain $J_{d_{1}+d_{2}-1}$. Since $J$ is strongly stable and $x_3$ is a strong Lefschetz element for $S/J$, they are
$x_{1}^{d_{1}+d_{2}-d_{3}-1}x_{2}^{2d_{3}-d_{1}-d_{2}}x_{3}^{d_{1}+d_{2}-d_{3}}$, $\ldots, 
x_{2}^{d_{3}-1}x_{3}^{d_{1}+d_{2}-d_{3}}$, therefore $J_{d_{1}+d_{2}-1}$ is the set
\[ \{ \{x_{1},x_{2}\}^{d_{1}+d_{2}-1},\ldots,x_{3}^{d_{1}+d_{2}-d_{3}}\{x_{1},x_{2}\}^{d_{3}-1},
x_{3}^{d_{1}+d_{2}-d_{3}+1}x_{1}^{d_{1}+d_{2}-d_{3}+1}\{x_{1},x_{2}\}^{2d_{3}-d_{1}-d_{2}-3},\] 
\[\ldots,x_{3}^{d_{1}}x_{1}^{d_{1}} \{x_{1},x_{2}\}^{d_{2}-d_{1}-1}, 
x_{1}^{d_{1}}x_{3}^{d_{1}+1} \{x_{1},x_{2}\}^{d_{2}-d_{1}-2},\ldots,x_{1}^{d_{1}}x_{3}^{d_{2}-1}\}.\]

One can easily prove by induction on $1\leq j \leq d_{3}-d_{2}$ that $J_{d_{1}+d_{2}-1+j}$ is the set
\[ Shad(J_{d_{1}+d_{2}-2+j})\cup \{
 x_{1}^{d_{1}+d_{2}-d_{3}+j-1}x_{2}^{2d_{3}-d_{1}-d_{2}-2j}x_{3}^{d_{1}+d_{2}-d_{3}+2j}, \ldots,
  x_{2}^{d_{3}-1-j}x_{3}^{d_{1}+d_{2}-d_{3}+2j}\} \]
Indeed, the case $j=1$ was already done and the induction step is similar. In particular, we get
\[J_{d_{1}+d_{3}-2} = \{\{x_{1},x_{2}\}^{d_{1}+d_{3}-2}, x_{3}\{x_{1},x_{2}\}^{d_{1}+d_{3}-3},\ldots,
                         x_{3}^{d_{1}+d_{3}-d_{2}-2}\{x_{1},x_{2}\}^{d_{2}},\]
\[x_{1}^{d_{1}}x_{3}^{d_{1}+d_{3}-d_{2}-1}\{x_{1},x_{2}\}^{d_{2}-d_{1}-1}, \ldots, x_{1}^{d_{1}}x_{3}^{d_{3}-2}\}.\]
Since $|J_{d_{1}+d_{3}-1}|-|Shad(J_{d_{1}+d_{3}-2})|=d_{1}$ we must add $d_1$ generators to $Shad(J_{d_{1}+d_{3}-2})$
in order to obtain $J_{d_{1}+d_{3}-1}$. Since $J$ is strongly stable and $x_3$ is a strong Lefschetz element for $S/J$, these new generators are 
 $x_{1}^{d_{1}-1}x_{3}^{d_{1}+d_{3}-d_{2}}x_{2}^{d_{2}-d_{1}},\ldots,x_{3}^{d_{1}+d_{3}-d_{2}}x_{2}^{d_{2}-1}$,so 
 $J_{d_{1}+d_{3}-1}$ is the set
\[\{ \{x_{1},x_{2}\}^{d_{1}+d_{3}-1},...,x_{3}^{d_{1}+d_{3}-d_{2}}\{x_{1},x_{2}\}^{d_{2}-1},
x_{1}^{d_{1}}x_{3}^{d_{1}+d_{3}-d_{2}+1}\{x_{1},x_{2}\}^{d_{2}-d_{1}-2},...,x_{1}^{d_{1}}x_{3}^{d_{3}-1}\}.\]

We prove by induction on $1\leq j\leq d_2-d_1$ that $J_{d_{1}+d_{3}-2+j}$ is the set
\[  Shad(J_{d_{1}+d_{3}-3+j})\cup\{ x_{1}^{d_{1}-1}x_{3}^{d_{1}+d_{3}-d_{2}+2j-2}x_{2}^{d_{2}-d_{1}-j+1} , \ldots, x_{3}^{d_{1}+d_{3}-d_{2}+2j-2}x_{2}^{d_{2}-j}\}. \]
Indeed, we already proved this for $j=1$ and the induction step is similar. We get
\[J_{d_{2}+d_{3}-2} = \{ \{x_{1},x_{2}\}^{d_{1}+d_{3}-2},x_{3}\{x_{1},x_{2}\}^{d_{1}+d_{3}-3},\ldots,
x_{3}^{d_{3}+d_{2}-d_{1}-2}\{x_{1},x_{2}\}^{d_{1}} \}.\] 
One can easily prove by induction on $1\leq j\leq d_1$ that 
\[J_{d_{2}+d_{3}-2+j} = Shad(J_{d_{2}+d_{3}-3+j})\cup \{
x_{3}^{d_{2}+d_{3}-d_{1}-2+2j}\{x_{1},x_{2}\}^{d_{1}-j}\} =\]\[ =  \{ \{x_{1},x_{2}\}^{d_{2}+d_{3}-2+j},\ldots, x_{3}^{d_{2}+d_{3}-d_{1}-2+2j}\{x_{1},x_{2}\}^{d_{1}-j}\}.\]
Finally, we obtain $J_{d_{1}+d_2+d_3-2}=S_{d_{1}+d_2+d_3-2}$ and therefore we cannot add new minimal generators of $J$ in degrees $>d_{1}+d_2+d_3-2$.
\end{proof}

\begin{cor}
In the above conditions of the above proposition, the number of minimal generators of $J$ is $d_1(d_2+1) - 
\left(\frac{\alpha}{2} \right)^{2} + 1$ if $\alpha$ is even or $d_1(d_2+1)-\frac{\alpha^{2}-1}{4} + 1$ if $\alpha$ is odd.
\end{cor}

\begin{ex} 

\begin{enumerate}
	\item Let $d_{1}=3$ , $d_{2}=5$ and $d_{3}=6$. Then
$ J = (x_1^{3}, x_1^{2}x_2^{3}, x_1x_2^{5},x_2^{6},$ 

$\;x_2^{4}x_3^{2}\{x_1,x_2\},\;x_2^{2}x_3^{4}\{x_1,x_2\}^{2},\; x_2x_3^{6}\{x_1,x_2\}^{2},\;x_3^{8}\{x_1,x_2\}^{2},\; x_3^{10}\{x_1,x_2\},\; x_3^{12})$. \vspace{1 pt}

  \item Let $d_{1}=4$ , $d_{2}=5$ and $d_{3}=6$. Then $ J = (x_1^{4},\;x_1^{3}x_2^{2},\;x_1^{2}x_2^{4},\;x_1x_2^{5},\;x_2^{7},\;x_3x_2^{6},$ 
  
  $\;x_2^{3}x_3^{3}\{x_1,x_2\}^{2},
\; x_2x_3^{5}\{x_1,x_2\}^{3},\; x_3^{7}\{x_1,x_2\}^{3},\; x_3^{9}\{x_1,x_2\}^{2},\; x_3^{11}\{x_1,x_2\},\;x_3^{13})$.
\end{enumerate}
\end{ex}

\begin{rem}
If $f_1,f_2,f_3\in S=K[x_1,x_2,x_3]$ is a regular sequence of homogeneous polynomials of given degrees $d_1,d_2,d_3$ such that $S/(f_1,f_2,f_3)$ has $(SLP)$, then the number of minimal generators of $J=Gin((f_1,f_2,f_3))$, $\mu(J) \leq d_1(d_2+1)+1$. This follows immediately from $2.4$, $2.9$, $3.4$, $3.9$, $3.14$ and $3.19$.
\end{rem}

\end{document}